\documentclass[11pt, A4page]{article}
\usepackage{amsmath}
\usepackage{graphicx}
\usepackage{comment}
\usepackage{amssymb}
\usepackage{pifont}
\usepackage{tabularx}
\usepackage{authblk}
\usepackage{natbib}
\usepackage[margin=1in]{geometry}
\usepackage[utf8]{inputenc}
\usepackage[english]{babel}
\usepackage{tikz}
\usepackage{tabularx}
\usepackage[flushleft]{threeparttable}
\usepackage{rotating}
\usepackage{multirow}
\usepackage{setspace}

\usepackage{imakeidx}
\makeindex[columns=3, title=Alphabetical Index, intoc]
\usepackage[titletoc]{appendix}
\usepackage{array}
\newcolumntype{P}[1]{>{\centering\arraybackslash}p{#1}}
\providecommand{\keywords}[1]{\textbf{\text{Keywords:}} }
\usepackage{authblk}
\usepackage{color}
\usepackage{mathtools}
\usepackage[utf8]{inputenc}
\usepackage[english]{babel}
\usepackage[inoutnumbered,linesnumbered,algoruled,slide,vlined]{algorithm2e}
\newtheorem{theorem}{Theorem}

\usepackage{ulem}
\usepackage{paralist} 

\usepackage{array}
\newcolumntype{P}[1]{>{\centering\arraybackslash}p{#1}}
\usepackage{float}

\begin{document}

\title
{
Dynamic Markdown Strategies for Perishable Products: \\ Waste Reduction versus Profitability}
\date{}
\author{
Nishika Bhatia\thanks{Jindal School of Banking \& Finance, O.P. Jindal Global University, India. Email: nbhatia@jgu.edu.in}
,\hspace{0.25cm}
Nalan G{\"{u}}lp{\i}nar\thanks{Durham University Business School, Durham, DH1 3LB, UK. Email: nalan.gulpinar@durham.ac.uk }
,\hspace{0.25cm}
Nursen Ayd{\i}n\thanks{Corresponding Author: The University of Warwick, Warwick Business School, UK. Email: naydin@wbs.ac.uk}
}
\vspace{-0.45cm}

\maketitle

\vspace{-0.5cm}

\begin{abstract}
{\bf Problem Definition:} 
 Freshness is a key attribute of perishable goods, particularly in the food sector, and strongly influences consumer purchasing behavior. Retailers often prioritize high freshness to remain competitive, yet maintaining this standard under demand uncertainty increases the likelihood of unsold inventory and waste.
Retailers use markdowns to mitigate this risk, but the timing and depth of discounts must be carefully coordinated to avoid lost revenue or excess waste. This paper studies how retailers can jointly optimize order and markdown pricing decisions to balance profitability and waste reduction for perishable products within an integrated decision-making framework. 

{\bf Methodology/Results:} The joint ordering-markdown pricing problem is formulated as a stochastic dynamic programming model that accounts for the demand cannibalization between regular price and markdown sales. The conflicting objectives of profit maximization and waste reduction are captured through a weighted objective function that reflects the retailer’s trade‑off preferences. To manage the complexity of the model, we propose an efficient solution approach and conduct numerical experiments to assess the performance of alternative markdown policies. 
Our results show that implementing multiple early price reductions with successive price reductions can simultaneously increase revenue and reduce waste.

{\bf Managerial Implications:} Our findings generate several key managerial insights: (i) firms can tailor their markdown policies according to their strategic priorities of profit maximization, waste reduction, or a balanced combination of both, (ii) early and well-timed markdowns based on available inventory levels allow firms to simultaneously attract quality-seeking customers and price-sensitive segments, (iii) implementing a second or multiple markdowns can further reduce waste while also increasing revenue.
\end{abstract}


\section{Introduction}

Perishable products (such as fruits, vegetables, dairy products, and cut flowers) are central to grocery retailing, accounting for more than half of total supermarket sales \citep{FMI2019,rts2022}. A recent report by \cite{NielsenIQ24} highlights that growth in perishable categories continues to outpace most other retail segments, underscoring their strategic importance for supermarkets. Consumers are highly sensitive to freshness, and empirical studies show that they assign measurable willingness-to-pay for even small extensions in shelf life, particularly in categories such as dairy and ready-to-eat meals \citep{Sousa25}. Furthermore, preferences on expiration dates strongly influence their purchase decisions \citep{Hansen24}. Consequently, retailers invest heavily in promoting the freshness, quality, and variety of their perishable offerings through targeted marketing and in-store displays \citep{Tescoo1}.

Although promoting perishables helps retailers establish a strong market position, it also introduces significant operational and environmental challenges. Displaying a wide variety of fresh-looking perishable products, including food items, often leads to excess inventory and, ultimately, food waste. 
According to the World Food Program USA \citep{WFP}, roughly one‑third of all food produced globally for human consumption (about 1.3 billion tons valued at around USD 1 trillion each year) is lost or wasted.
Recent global estimates from the UN Environment Program indicate that a considerable proportion of total food waste arises in the downstream stages of the supply chain, with approximately 12\% occurring at the retail level \citep{UNEP24}. In the United States alone, around 30\% of the fresh products available from supermarkets are discarded each year, resulting in roughly 16 billion pounds of food waste; the financial loss from these wasted products is estimated to be nearly twice the profit generated from total food sales \citep{rts2022}. This imbalance between freshness and profitability has turned perishable inventory management into a crucial sustainability and operational challenge.

Recognizing the scale of this problem, the United Nations has committed the global community to halve per capita food waste at the retail and consumer levels and reduce food loss across supply chains by 2030 \citep{UNEP24}. Achieving this ambitious target requires action by policymakers, supply chain actors, and retailers. To minimize wastage and reduce potential losses, many retailers have begun to implement price reductions on products near the end of their shelf life. Traditionally, markdown strategies typically involve a single (one-time) price reduction applied shortly before the expiration date of the product \citep{Chehrazi25}. However, a recent study by McKinsey emphasizes a shift towards earlier and more frequent price adjustments throughout the product lifecycle \citep{Bahr22}. For instance, Hoogvliet supermarkets in the Netherlands have partnered with Wasteless to implement a dynamic markdown system that gradually lowers prices over time rather than applying a single discount on the last day. This system begins with small markdowns (e.g., 2–5\%) and gradually increases them (e.g., 10–15\%) closer to expiration \citep{Blum21, Morrison22, RetailCare25, trellis24}. Similar practices have been documented in different retail settings. The S-market in Finland employs a two-tier markdown system: short-dated products such as meat, fish, and vegetables are reduced by 30\% during the day, followed by a 60\% discount after 9 pm to encourage end-of-day clearance \citep{WEF2019}. \cite{Hansen24} conducted a large-scale field experiment with a European grocery retailer that implemented a dynamic pricing based expiration date. The system displayed both the remaining shelf life and the corresponding discount for the older products, allowing consumers to compare use-by dates and prices in real time. In their experiment, the retailer placed the products near expiration on the front of the shelf to encourage the selection by price-sensitive consumers. Their findings show that combining shelf rotation with dynamic markdowns can reduce retail food waste by up to 30\%. These practices demonstrate that age-differentiated, multi-tiered markdowns are now becoming mainstream through the integration of real-time inventory and pricing systems. 
In addition, emerging technologies such as electronic shelf labels, pricing software, and mobile apps are increasingly enabling and supporting the use of age‑based dynamic markdown pricing.

Markdown sales, however, inherently involve a tension between two objectives: maximizing revenue and minimizing waste. They can substantially reduce waste from unsold inventory, but may also cannibalize demand for full-priced products. Setting a very low price for old inventory can result in a revenue loss, as consumers who initially intend to buy new inventory might switch to older inventory due to an attractive price difference. In contrast, pricing old inventory too high may decrease demand; consequently, this leads to higher waste costs. The central managerial challenge, therefore, lies in determining the timing and intensity of the markdown that best balance revenue generation and waste reduction under demand uncertainty. Several studies, such as \citep{Bahr22} and \cite{retalon2022}, specifically emphasize the importance of developing joint markdown and replenishment strategies to improve profits. Coordinating markdown strategies with replenishment decisions is therefore critical to minimize the economic burden from mismatches between perishable supply and demand.

Motivated by these practical developments, we consider the joint decision-making problem of a retailer selling a perishable product with a finite shelf life within a periodic-review framework. Specifically, we focus on how integrated ordering and pricing decisions can help balance profitability and waste reduction. The replenishment cycle is shorter than the shelf life of the product; therefore, new (fresh) and older inventories coexist on the shelves simultaneously. We define the age of a product as the number of periods that have passed since its arrival at the store, and any unit that remains unsold until its expiration is discarded, incurring a disposal cost. The retailer can conduct markdown sales for older inventory to reduce waste, always pricing aging products below fresh ones. The demand in each period is uncertain and any unmet demand due to stock shortages is lost. Customers are heterogeneous in their preferences for freshness and price. Although some prefer to purchase fresh products at the regular price, others opt for older discounted products, choosing the option that provides the highest utility. The impact of demand cannibalization between regular and markdown inventories is captured by evaluating dynamic changes in customer choices over time.  The retailer determines both the quantity of fresh products to order and the markdown prices to charge for the aging inventories in each time period. 

We formulate this problem as a stochastic dynamic programming model. The inherent challenge in managing perishable inventory lies in the dimensional expansion of the state space, as a result of the continuous tracking of the aging inventories. Due to the complexity of the dynamic decision-making problem, existing research focuses predominantly on inventory systems with a two-period lifetime product. Recent discussions in the operations management literature emphasize that as detailed inventory age data become more accessible, there is an increasing need for models and algorithms that can explicitly link product age, markdown timing, and demand cannibalization effects \citep{akkacs2022reducing}. Our research  responds directly to this call by developing a general framework for perishable products that integrates these dimensions. 
This study advances the existing literature in several important directions. First, the proposed stochastic dynamic programming model accommodates products with general $n$-period lifetimes, extending beyond the two-period settings that dominate prior research. It explicitly tracks the evolution of inventory by product age and allows multiple (gradual) markdowns based on remaining shelf life, consistent with recent developments in retail markdown and inventory practices. This formulation enables the systematic analysis of alternative pricing and ordering strategies, including the timing and magnitude of markdowns across different inventory ages. However, tracking inventory by product age over time creates a high-dimensional state space, making the problem analytically challenging. To manage this complexity, we exploit key structural properties of the dynamic programming model. We show that the value function is $k$-concave and supermodular in inventory levels and leverage these properties to design an efficient solution algorithm that yields optimal ordering and markdown policies. Second, our model incorporates customer choice and demand cannibalization, capturing heterogeneity in consumer sensitivity to freshness and price. This choice-based representation links markdown decisions to consumers’ purchasing behavior, allowing the model to evaluate how price adjustments across product ages shape overall demand and waste levels. Finally, through extensive numerical experiments, we examine how early and multiple markdown strategies influence the trade-off between profitability and waste reduction. Contrary to the common perception among grocery retailers that earlier or more frequent discounts inevitably erode revenue, our results reveal that gradual markdowns can preserve profitability while reducing waste.

The remainder of the paper is organized as follows. Section 2 focuses on reviewing the literature by providing details of existing studies relevant to our research. The stochastic dynamic programming formulation of the joint ordering-markdown pricing problem is presented in Section 3. In Section 4, we introduce an efficient solution approach that takes advantage of  specific features of the underlying stochastic dynamic program. We show that this approach provides an optimal ordering-markdown pricing strategy. 
Section 5 focuses on the design of numerical experiments and presents computational results that help to derive managerial insights. The concluding remarks are provided in Section 6. 

\section{A Brief Literature Review}

Research on perishable inventory management has attracted significant attention over the years. Several comprehensive reviews synthesize progress in this area, including \cite{nahmias11}, \cite{karaesmen2011managing}, \cite{bakker2012review}, \cite{Janssen16}, and \cite{akkacs2022reducing}. These studies review fundamental models for managing inventories of items with limited shelf life, addressing challenges such as the trade-off between stock availability and profitability, demand management, shelf-space constraints, and spoilage. Building on this foundation, our study contributes to the growing stream of research integrating inventory management and dynamic pricing for perishable products, with particular emphasis on how pricing and ordering decisions can affect product wastage. In this area, the main focus has been towards the periodic review models where inventory has been tracked at specific periods rather than continuously. The reader is referred to 
\cite{chen2021optimal}, \cite{li2022joint} and \cite{vahdani2022coordinated} for different continuous review models introduced for perishable products with random lifetime. 

We organize our review of the literature around three key themes and explicitly highlight the evolution of the research focus in this area. First, we examine the literature on joint inventory management and pricing models for perishable products. \cite{ferguson2007should} analyze inventory management and pricing decisions for a company that sells a food product with an exact two-period life cycle. The procurement and pricing decisions for fresh products are made in the first time period considering demand uncertainty. At the start of the second time period, the decisions are to determine how much leftover inventory of the old product to carry over, how much fresh products to procure, and what prices to charge for both new and old inventories by assuming that demand in the second period is deterministic. \cite{li2009note} study joint pricing and inventory control for a two-period lifetime perishable product over an infinite horizon. They analyze the structural properties of the optimal inventory replenishment and pricing policy by assigning a single price for both old and new inventories. In this model, although inventories that are tracked differently age, there is no price differentiation between old and new inventories. Because exact solutions become complex for longer lifetimes, they propose a heuristic that fixes both a replenishment level and a selling price. Price differentiation for perishable products with a general shelf life of $n$ periods is addressed by \cite{li2012grocery}, where in each period the retailer makes replenishment and pricing decisions and decides whether to retain or discard the remaining stock. In this model, old and new inventories are never sold simultaneously, and the old inventory is discarded whenever a new order is placed. \cite{chen2013joint} extend the work of \cite{li2009note} by considering various inventory consumption scenarios, specifically examining the implications of adopting the first-in-first-out and last-in-first-out approaches. Through their analysis, they obtain structural insights into the optimal inventory management policy. In a related work, \cite{chen2014coordinating} investigate how the lead time affects the replenishment decisions for perishable products. 
According to the existing literature, old and new inventories carry the same price, but any excess inventory is disposed of prior to expiration.
\cite{herbon2017should} examines whether retailers benefit from allowing old and new perishable inventory to coexist on shelves at different price points within a deterministic demand framework. Based on this, \cite{atan2023displaying} has recently analyzed the possible impact of demand uncertainty and determined the optimal strategies for the shelf presentation of perishable goods with a two-period shelf life. 

Several studies have specifically examined inventory policies that account for product age, even if prices are not the primary focus. \cite{deniz2010managing} study inventory issuing and replenishment policies for perishable products with a two-period lifespan. They investigate how demand might depend on the remaining shelf life of a product and evaluate different issuing policies. \cite{li2014multi} study the structural properties of optimal inventory policies for perishable products where the retailer sells both newly stocked and aged products at differentiated prices. At each time period, the retailer determines the quantity of new products to order, the volume of older unsold inventory to retain, and the amount of product to offer at the clearance sale given the regular and clearance sale prices. They assume that clearance products will be sold out due to their high demand, leading the firm to strategically decide the amount of older products to offer at clearance sales. \cite{Zhang20} and \cite{chen2021} extend this line of research considering similar periodic-review settings in which the retailer observes the age distribution of inventory, determines replenishment quantities, and decides how much aged stock to clear or dispose of before expiration. Because the resulting dynamic programs are computationally intractable, both studies develop tractable approximation approaches for ordering and disposal decisions, achieving near-optimal performance. Complementing these analytical models, \cite{Keskin22} develop a data-driven dynamic pricing and ordering framework for perishable inventory operating in a changing demand environment. They assume that the demand–price relationship is unknown and may vary over time. The prices of new and old inventory are kept identical and the retailer learns optimal pricing and ordering policies through a dynamic data-driven program. These studies lay the groundwork for understanding how inventory dynamics and demand uncertainty interact in perishable contexts, but they do not yet fully capture the strategic role of markdown pricing and its interaction with inventory, a gap that more recent research is beginning to address.

The second stream of literature examines how products of different ages coexist and compete in the market, introducing markdown pricing and intertemporal cannibalization mechanisms. \cite{sainathan2013pricing} studies a joint pricing and inventory model of a perishable product with a two-period shelf life in which customers choose between fresh and older items based on linear utilities, capturing the trade-off between freshness and price. Their results show that selling older products can become profitable under uncertainty, especially when demand variability is high. \cite{hu2015joint} analyze strategic customer behavior, where consumers may delay purchases in anticipation of discounts, creating intertemporal cannibalization. However, their model is distinctly designed for products with a one-period lifespan (such as bakery goods and other food items), addressing the unique challenges associated with managing highly perishable inventory. \cite{li2016managing}, in contrast, address the problem of joint replenishment and clearance sales problem for products with a general shelf life of $n$-period, highlighting threshold-based clearance policies and the role of the oldest inventory. 
Both \cite{hu2015joint} and \cite{li2016managing} assume exogenously given static prices for old and new inventories. Based on this line of work, \cite{chua2017optimal} and \cite{fan2020dynamic} develop inventory and pricing strategies for a two-period perishable product. The research of \cite{chua2017optimal} focuses on determining the discounts for older merchandise, while \cite{fan2020dynamic} explore pricing strategies for both fresh and aged products. In these studies, pricing decisions and customer demand for the perishable product are assumed to be independent of time. 
\cite{chao2015approximation} address the time-dependent demand for products with a $n$-period lifespan, but their emphasis lies in developing approximation algorithms for tractable inventory policies rather than pricing strategies. \cite{chintapalli2015simultaneous} develops a joint inventory-pricing decision model for a retailer facing substitutable demand for multi-period food products, assuming that consumers can choose between old or new products. Due to the complexity of the underlying problem, the model is simplified to consider a scenario in which the discounted prices for older food items remain constant over time.  \cite{Moshtagh25} study a dynamic inventory and pricing control problem for a perishable product that progresses through multiple shelf-life phases. They first formulate a continuous-review stochastic control model to describe the joint evolution of inventory and pricing decisions over time. To obtain a tractable representation, the authors discretize the model into a finite-horizon dynamic program in which each period is short, and in each state, the firm decides on the prices for all freshness levels and whether to produce. To address the computational complexity of the model, particularly in large-scale systems, the authors establish structural properties of the optimal policy and propose several heuristics that deliver near-optimal performance. 

The third stream of literature focuses on waste management and sustainability in perishable inventory systems, emphasizing how operational and pricing decisions can mitigate product spoilage. Although most previous studies consider spoilage at the end of lifetime of a product, explicit attention to waste management has been largely overlooked. Most articles assume that the unsold stock is discarded without penalty or has a constant salvage value \citep{abouee2019, Zhang20, chen2021}. When designed appropriately, markdowns can accelerate the sale of aging inventory that would otherwise expire, creating a potential win–win outcome that improves both profitability and sustainability. \cite{Boer22} study markdown pricing and inventory control for perishable products with a two-period shelf life, explicitly treating wastage in the objective function along with profit maximization. They consider a retailer that restocks daily to a certain level and faces customers whose purchase probabilities depend on both the price and the remaining shelf life of the product. By analyzing the steady-state of this system, \cite{Boer22} show that a markdown policy on the last day of shelf life of products produces strictly less waste than a fixed price policy in a broad set of scenarios. \cite{Zhou23} develop a two-period dynamic pricing and inventory model for a perishable product with two replenishment periods designed to satisfy demand throughout the entire sales cycle. A distinctive feature of the study is its focus on sustainability objectives, as the model jointly optimizes ordering, transfer, and return decisions throughout the two periods to rationalize replenishment and reduce overproduction. The authors demonstrate that coordinated inventory-adjustment policies can simultaneously decrease surplus and shortage levels, thus mitigating product waste while improving overall efficiency. \cite{Sanders24} extends this stream by analyzing the relative effectiveness of dynamic pricing and organic waste bans using an empirical dynamic newsvendor model calibrated with retail grocery data. The study shows that dynamic pricing can reduce waste by more than 20\% while increasing retailer profits, whereas waste bans achieve only modest reductions at the expense of profitability and welfare. \cite{Chehrazi25} analyze transaction-level data from a major U.S. grocery retailer chain and reveal that prices for perishable goods often remain inflexible as expiration approaches. 
This is mainly attributed to informational frictions between pricing and inventory systems. Their findings indicate that advancing markdowns even by one day can lead to waste reduction, while maintaining or slightly enhancing revenue performance. Complementing these studies, \cite{akkacs2022reducing} highlight that as the field continues to evolve, there is growing interest in models that integrate aging inventories, markdown timing, and cannibalization effects to improve perishable inventory management. In particular, they emphasize that traceability systems can identify items approaching expiration, enabling targeted promotions and improved stock rotation, thereby reducing waste. Our study responds directly to this call by developing a stochastic dynamic framework that jointly optimizes ordering and markdown decisions while accounting for demand cannibalization under uncertainty.
\vspace{-0.95cm}
\begin{center}
\begin{threeparttable}[htp!]
\caption{\label{tab:2-ch3} \text{Classification of the most relevant research papers}}
\scriptsize
\begin{tabular}{p{3.2cm}|cccc|ccc|cccc|cc|c|}
\cline{2-15} 
\multicolumn{1}{c|}{}& 
\multicolumn{4}{|c|}{\text{\textbf{ Decisions }}} & \multicolumn{3}{c|}{\text{\textbf{ Parameters }}} &  \multicolumn{4}{c|}{\text{\textbf{ Demand Model }}} &  \multicolumn{2}{c|}{\text{\textbf{ Solution }}} &  \multicolumn{1}{c|}{\text{\textbf{Waste}}}
\\ 
\cline{2-15} 
\cline{2-15} 
\text{Research Papers }  
& \begin{turn}{90} 
{\it Ordering (S, D)}
\end{turn}& \begin{turn}{90} 
{\it Pricing New Product (S, D)} 
\end{turn} & \begin{turn}{90} 
{\it Pricing Old Product (S, D)} 
\end{turn}& \begin{turn}{90} 
{\it Markdown  (F, T, D)}  \;
\end{turn}
& \begin{turn}{90} 
{\it Multi-age Inventory}
\end{turn} 
& \begin{turn}{90} 
{\it Age (1, 2, n)}
\end{turn}
& \begin{turn}{90} 
{\it Planning Horizon}
\end{turn}
& \begin{turn}{90} 
{\it Demand (U, C)}
\end{turn}
& \begin{turn}{90} 
{\it Price Independent}
\end{turn} & \begin{turn}{90} 
{\it Linear/Log/Iso-elastic} \;
\end{turn} & \begin{turn}{90} 
{\it Customer Utility}
\end{turn}
& \begin{turn}{90} 
{\it Optimal}
\end{turn} & \begin{turn}{90} 
{\it Approximation}
\end{turn}
& \begin{turn}{90} 
{\it Waste-in-Objective (P)} \;
\end{turn}\\
\hline
\hline
\textsl{Ferguson et al. (2007)} & \textsl{D} & \textsl{D} & \textsl{S} & \textsl{F(s)} & $\bullet$ & 2 & 2 & {C} & $\bullet$ & -- & $\bullet$ & $\bullet$ & -- & -- \\
\textsl{\cite{li2012grocery}} & \textsl{S} & \textsl{D} & \textsl{D}  & -- & -- & $n$ & $\infty$ & \textsl{U} & $\bullet$ & $\bullet$ & -- & -- & $\bullet$ & \textsl{P} \\
\textsl{\cite{sainathan2013pricing}} & \textsl{D} & \textsl{D} & \textsl{D} & -- & $\bullet$ & 2 & $\tau$ & \textsl{U} & $\bullet$ & -- & $\bullet$ & -- & $\bullet$  & -- \\
\textsl{\cite{li2014multi}} & \textsl{D} & -- & -- & --& $\bullet$ & $n$ & $\tau$ & \textsl{U} & $\bullet$ &  $\bullet$ & -- & -- & $\bullet$ & \textsl{P}  \\
\textsl{\cite{chen2014coordinating}} & \textsl{D} & \textsl{D} & \textsl{D} & -- & $\bullet$ & $n$ & $\tau$ & \textsl{U} & $\bullet$ &  $\bullet$ & -- & -- & $\bullet$ & \textsl{P}\\
\textsl{\cite{hu2015joint}} & \textsl{D} &-- & -- &\textsl{T(s)} & -- & 1 & $\infty$ & \textsl{C}  & -- & -- & -- &  $\bullet$ & -- & -- \\
\textsl{\cite{chao2015approximation}} & \textsl{D} & -- & -- & -- & $\bullet$ & $n$ & $\tau$ & \textsl{U} &  --& --& --& --& $\bullet$ & \textsl{P}\\	
\textsl{\cite{li2016managing}} & \textsl{D} & --& --& \textsl{T(s)} & $\bullet$ & $n$ & $\tau$ & \textsl{U} & -- &-- & --& --& $\bullet$ & \textsl{P} \\	
\textsl{\cite{chua2017optimal}} & \textsl{D} & -- & \textsl{D} & \textsl{T(m)} & $\bullet$ & 2 & $\tau$ &  \textsl{U} & $\bullet$ & $\bullet$  & -- & --&  $\bullet$  & --  \\
\textsl{\cite{herbon2017should}} & \textsl{S} & \textsl{S} &  \textsl{S} & -- & $\bullet$ & 2 & $\tau$ & \textsl{C}  & -- & -- & -- & $\bullet$ &  -- & -- \\
\textsl{\cite{fan2020dynamic}} &  \textsl{S} & \textsl{D} & \textsl{D} & -- & $\bullet$ & 2 & $\tau$ & \textsl{U}  & $\bullet$ & -- &$\bullet$ & -- & $\bullet$ & -- \\
\textsl{\cite{Zhang20}} &  \textsl{D} & -- & -- & -- & $\bullet$ & $n$ & $\tau$ & \textsl{U}  & $\bullet$ & -- & -- & -- & $\bullet$ & \textsl{P} \\
\textsl{\cite{chen2021}} &  \textsl{D} & -- & -- & -- & $\bullet$ & $n$ & $\tau$ & \textsl{U}  & $\bullet$ & -- &-- & -- & $\bullet$ &\textsl{P} \\
\textsl{\cite{Boer22}} & \textsl{S} & \textsl{S} & \textsl{S} & \textsl{F(s)} & $\bullet$ & 2 & $\infty$ & \textsl{U}  & -- & -- &$\bullet$ & -- & $\bullet$ & \textsl{P} \\
\textsl{\cite{Zhou23}} & \textsl{D} & \textsl{D} & \textsl{S} & \textsl{F(s)} & $\bullet$ & 2 & 2 & \textsl{U}  & $\bullet$ & -- &$\bullet$ & -- & $\bullet$ & \textsl{P} \\
\textsl{\cite{Moshtagh25}} & \textsl{S} & \textsl{D} & \textsl{D} & \textsl{D(m)} & $\bullet$ & n & $\tau$ & \textsl{U}  & -- & -- &$\bullet$ & -- & $\bullet$ & \textsl{P} \\
\hline  \hline
\textsl{This Paper} & \textsl{D} & \textsl{D} & \textsl{D} &  \textsl{D(m)}  & $\bullet$ & $n$ & $\tau$ & \textsl{U}  & -- & -- &$\bullet$ & $\bullet$ & -- & \textsl{P} \\
\hline
\end{tabular}
\begin{tablenotes}
{\item S: Static, D: Dynamic, F: Fixed, T: Threshold, C: Certain, U: Uncertain, P: Penalty, s: Single, m: Multiple, --: Nonuse}
\end{tablenotes}
\end{threeparttable}
\end{center}
 To provide a structured overview of the studies discussed so far, Table \ref{tab:2-ch3} summarizes the key characteristics of the most relevant research on joint inventory management and pricing for perishable products. The table classifies the studies according to the modeling scope and solution structure, reflecting how inventory aging, demand uncertainty, computational tractability, and waste management are addressed in existing work. The first four columns in ``Decisions" indicate whether the ordering, pricing, and markdown policies are modeled as static (S), dynamic (D), fixed (F) or threshold-based (T). The notation \textit{s} and \textit{m} in the column ``Markdown" indicates whether the study considers a single markdown or multiple markdowns on the same product. The next three columns under ``Parameters" summarize how each study characterizes multi-age inventories, the product lifetime (one, two, or $n$ periods), and the planning horizon (finite or infinite). The following block titled ``Demand Model" outlines how demand is modeled in each study, specifying whether it is uncertain or certain, price-dependent, follows a particular functional form (linear, logarithmic, or iso-elastic), or incorporates customer utility to capture substitution or freshness preferences. The block ``Solution"  distinguishes between studies providing optimal analytical solutions and those relying on approximation or heuristic approaches to overcome model complexity. Finally, the last column indicates whether the waste cost is explicitly included in the objective function (as a penalty) or treated implicitly.

In general, our research differs from the existing literature in several important ways. First, most previous studies develop joint pricing and inventory policies within the restrictive setting of products with a shelf life of two periods, focusing on static or simplified markdown structures \citep{ferguson2007should, li2009note, herbon2017should,  fan2020dynamic}. Second, while a few papers extend these frameworks to multi-period lifetimes, they typically overlook the explicit role of inventory freshness in markdown pricing or treat markdown prices as exogenously fixed \citep{chintapalli2015simultaneous, li2016managing, chua2017optimal}. By contrast, our study explicitly associates different pricing decisions with old and new inventories for a general $n$-period lifetime product and models the markdown process dynamically over time under demand uncertainty. Third, although a growing stream of work has begun to incorporate sustainability considerations and waste-reduction objectives \citep{Boer22, Zhou23, Sanders24}, these models often address ordering or pricing decisions in isolation rather than as part of an integrated ordering–markdown framework. Building on the agenda proposed by \cite{akkacs2022reducing}, we contribute to this literature by developing a unified model that jointly optimizes ordering and markdown pricing decisions, accounts for demand cannibalization between fresh and old inventories, and provides tractable structural results that enable efficient solution algorithms. 

\section{ Dynamic Ordering and Markdown Pricing Model}

In this section, we present a formulation of the joint ordering and markdown pricing problem of a retailer (such as Tesco and Aldi) using a stochastic dynamic program (SDP) under a customer choice model. The notation used in the problem formulation is summarized in Table \ref{tab:var0}. The firm sells a perishable product over a finite planning horizon, indexed by $t = 1, \cdots, \tau$. 
 Our focus is on perishable products with a short lifetime (including milk, fruits, vegetables, and other grocery items) where it is common to observe different expiration dates simultaneously. The product has a lifetime of $n$ periods (i.e., it perishes after $n$ time periods). The firm regularly reviews the available inventory and places an order for fresh products. We assume that the new inventory (fresh products) can be displayed alongside the older inventory on the same shelves. Let $i$ denote the age of the product. Specifically, $i=0$ refers to a fresh product, while older inventories at various ages are represented by $i \in \{1, \cdots, n\}$. Demand for available inventories at time $t$ is assumed to be uncertain and denoted by a random variable $\Tilde{d}_{t}$ for $t=1, \cdots, \tau$. Customer demand at any time period is fulfilled from available inventories using a first-in first-out (FIFO) issuing policy. 
 In other words, the oldest inventory is sold first when a customer's request arrives. If customer demand is not met due to inventory unavailability, the firm incurs a shortage penalty of $\gamma$ per unit. In addition, unsold inventories are carried over to the next period by reducing the remaining lifetimes by one, and incur a holding cost of $h$ per unit. The expired inventories are disposed of with a unit waste cost of $\eta$.

\begin{table}[h]
\centering
\footnotesize{ \caption{ Description of notation}
\label{tab:var0}
\begin{tabular}{ll}
\hline
&\text{\it Model parameters}  \\ \hline
$\tau$ & planning horizon discretized by time periods $t = 1, \cdots, \tau$ \\ 
$ n$ & fixed age of perishable product \\
$\kappa $ & production capacity  \\ 
$\delta$  & discount factor \\ 
$r$ & unit regular price \\
$c$ &  unit order cost  \\ 	
$h$ & per-unit holding cost \\ 
$\eta$ & unit wastage cost for discarded product \\
{\color{black} $\gamma$} & {\color{black} unit shortage cost for unmet demand} \\
$\tilde{d}_{t}$ & uncertain demand at time $t$ \\
\hline
&\text{\it State variables} \\
\hline
$\mathbf{x_{t}}$ & vector of regular inventory levels $x_{i,t}$ of age $i$ at the beginning of time $t$  \\
$\mathbf{y_{t}} $  & vector of markdown inventory levels $y_{i,t}$ of age $i$  at the beginning of time $t$  \\
$\mathbf{I_{t}} $  & inventory vector consisting of regular and markdown inventory as $\mathbf{I_{t}} = (\mathbf{x_{t}}, \mathbf{y_{t}}) $	\\
$\mathbf{w_{t}}$  & vector of binary variables for markdown information $w_{i,t}$ for inventory age $i$ at time $t$  \\ \hline
&\text{\it Actions} \\
\hline
$p_{t} $  & markdown price  at time $t$  \\
$q_{t}$ & amount of product to order at time $t$ \\ \hline 
\end{tabular}
}
\end{table}

{\bf System Dynamics:} 
 The firm can implement markdown sales by reducing the prices of unsold and older inventory to reduce potential losses and minimize waste. This practice of adjusting prices based on the age of inventory creates internal competition, as discounted items can draw demand away from products sold at regular prices. We assume that the retailer may apply a reduction to inventories of any age at any time. Consequently, multiple markdowns can be applied to products of varying ages (so called age-differentiated products), requiring inventory and pricing to be tracked across different age classes over time.
This framework accommodates early markdowns prior to expiration and allows for the continuous introduction of fresh products onto shelves at any time. In addition, we assume that if inventory is marked down at a time, it cannot be marked up in future time periods. Thus, to keep track of the price at any given time, we categorize the inventory as {\it regular} or {\it markdown}. In addition, a two-price structure (as regular and markdown) is implemented to reflect customer preferences for fresh versus discounted items.

Let $x_{it}$ and $y_{it}$ denote the quantities of regular and markdown inventories of age $i$ at the beginning of time $t$, respectively. We define a state vector $\mathbf{I}_{t} = (\mathbf{x}_{t}, \mathbf{y}_{t})$ 
consisting of regular inventories $\mathbf{x}_{t} = \{x_{it} \in \mathbb{Z}^{+}: i=0,\cdots,n \}$ and markdown inventories $\mathbf{y}_{t} = \{y_{it} \in \mathbb{Z}^{+}: i=0,\cdots,n \}$. To keep track of whether inventory of a given age is currently marked down, we introduce a vector of binary variables $\mathbf{w_{t}} = \{w_{it} \in \{0,1\}: i=0,1,\cdots,n+1\}$ where $w_{it} = 1$ indicates that inventory at age-$i$ is sold at the markdown price in period $t$. A two‑price system, as regular and marked down, reflects the preference of customers for fresher products versus older discounted items. Thus, a state of the dynamic system at time $t$ consists of amounts of (regular and markdown) inventories $\mathbf{I}_{t}$ and markdown indicators $\mathbf{w}_{t}$. It should be noted that, because the decision model treats regular and markdown inventories as separate pathways, maintaining distinct state variables for each category is necessary to track their evolution accurately. Moreover, for any period $t$ and age $i$, at most one of  $x_{it}$ or $y_{it}$ can be positive; they cannot be positive simultaneously. This mutual exclusivity ensures that, even though the model tracks two inventory types, the effective state space remains compact. 

\begin{figure}[h!]
\centering
\includegraphics[scale=0.18]{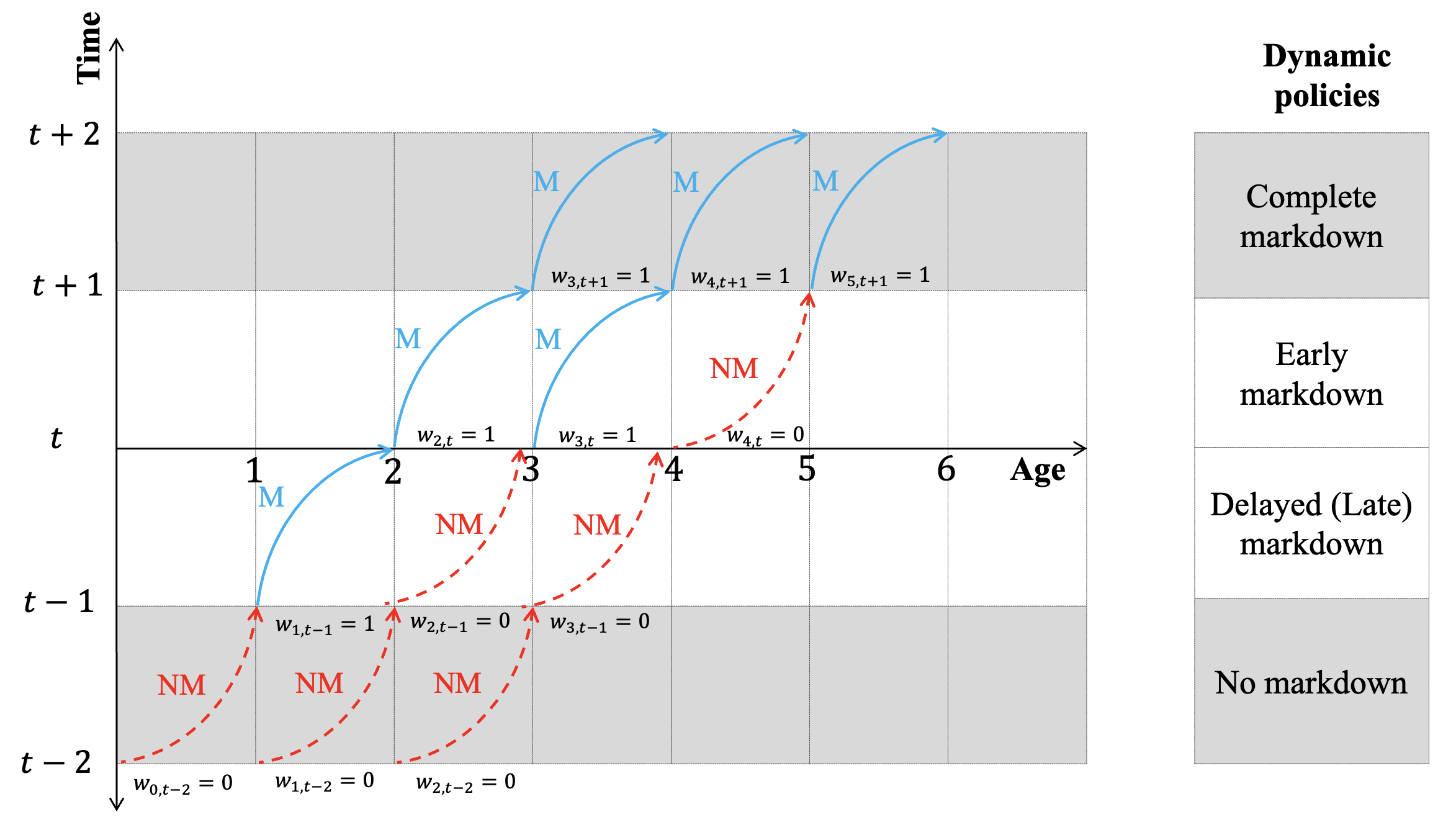}
\caption{State transition process involving markdown information over time}
\label{markdown}
\end{figure}
\vspace{-0.1cm}
We can now outline the state transition equations at each time. As shown in Figure \ref{markdown}, a drop in inventory state $w_{i,t+1}$ with age $i$ at time $t+1$ progresses from state $w_{i-1,t}$. Thus, the transition from $w_{i-1,t}$ to $w_{i,t+1}$ captures markdown information, indicating whether the remaining inventory aged $i-1$ is subject to price reduction or is sold at regular price at time $t$. For example, a transition from $w_{4,t} = 0$ to $w_{5,t+1} = 1$ means that inventory aged $4$, sold at regular price at time $t$ (no markdown - $NM$), becomes  marked down (M) at time $t+1$ and aged 5. Once the price of perishable inventories is reduced, they remain under markdown (i.e., never sold at its original price again) during their remaining shelf life. This implies that if the inventory at age $i-1$ is marked down $w_{i-1,t} = 1$ at time $t$, then it remains  marked down $w_{i,k} = 1$ for $k = t+1, t+2, \cdots, \tau$. On the other hand, if the inventory is sold at a regular price at time $t$, it can be offered at a regular or discounted price in the next period $t+1$. 
Thus, the markdown state transition for $i=1,\ldots,n+1$  can be expressed as follows;
\begin{equation}
\begin{aligned}
\displaystyle w_{i,t+1} = \left\{
\begin{array}{ll}
    1      &\text{ if } w_{i-1,t} = 1  \\
    0 \;\; \text{or} \;\; 1   &\text{ if } w_{i-1,t} = 0
\end{array}
\right\}.
\end{aligned}
\end{equation}
Throughout the planning horizon, various dynamic markdown strategies (namely 
{\it complete}, {\it partial}, and {\it no markdown} as illustrated in Figure 1) can be observed. Under the complete price reduction strategy, all inventory items, regardless of age, are offered at reduced prices. In contrast, the no-markdown policy maintains regular pricing for all inventory. The partial markdown strategy (implemented in early or delayed (late) states)  involves a mix in which some inventory is sold at reduced prices while others remain at regular price. For example, at time $t-1$  inventories with age-1 can be discounted (as $w_{1,t-1} = 1$), while inventory for age $i=2, \cdots,n$ is displayed at its regular price as ($w_{i,t-1} =0$). The quantities of regular and markdown inventories, denoted by $x_{i,t}$ and $y_{i,t}$, respectively, with age $i$ at time $t$ evolve from the remaining inventory $x_{i-1,t-1}$ and $y_{i-1,t-1}$ of age $i-1$ at previous time $t-1$. Assuming a first-in, first-out issuing policy, the transitions of regular and markdown inventories of age $i$ from $t-1$ to $t$ are governed by the following equations, respectively;

\setlength{\belowdisplayskip}{0pt} \setlength{\belowdisplayshortskip}{0pt}
\setlength{\abovedisplayskip}{0pt} \setlength{\abovedisplayshortskip}{0pt}
\begin{equation}\label{eqn:1-ch3}
\begin{aligned}
x_{i,t} &= max\Big\{ (1-w_{i,t})x_{i-1,t-1} - max \Big\{\tilde{d}_{t-1}  - \sum_{j=i}^{n-1} \left(1-w_{j+1,t}\right) x_{j,t-1}, \;0  \Big\}, \;0 \Big\}\; \text{and} \\
y_{i,t} &= max\Big\{ y_{i-1,t-1} + w_{i,t}\; x_{i-1,t-1} - max\Big\{\tilde{d}_{t-1}  - \sum_{j=i}^{n-1} (y_{j,t-1} + w_{j+1,t}\; x_{j,t-1}), \;0  \Big\}, \;0 \Big\}.
\end{aligned}
\end{equation}
Taking into account the markdown status of the inventories, the amount of regular inventory aged $i-1$ at time $t-1$ is calculated as $(1-w_{i,t})\; x_{i-1,t-1}$. On the other hand, the quantity of markdown inventory that ages $i-1$ at time $t-1$ is formulated as $y_{i-1,t-1}+w_{i,t} \; x_{i-1,t-1}$.

The total quantity of inventory offered at regular or markdown prices, respectively, independent of age categories, can be formulated as
\begin{equation}\label{eqn:1-ch3}
\begin{aligned}
X_{t} =  \sum_{i=0}^{n} (1-w_{i+1,t+1})x_{it}, & \; \; \text{and} \; &
Y_{t} = \sum_{i=0}^{n} y_{it} + w_{i+1,t+1}x_{it}. 
\end{aligned}
\end{equation}

{\bf Actions:} 
At the beginning of each time, the company reviews its  inventory in different ages and can order fresh products. Simultaneously, they determine markdown prices for the aging inventory. We assume zero lead time for ordering, meaning that fresh products are received instantaneously. Let $\kappa$ represent the fixed ordering capacity. 
{\color{black}
We incorporate a capacitated order quantity in our model to reflect the practical constraints faced by retailers, such as limited shelf space and warehouse capacity.
To ensure that the quantity ordered at time $t$ does not exceed this capacity, we impose the condition}
\begin{equation} \label{con1}
0 \leq  q_t \leq \kappa. 
\end{equation}
Inventory can be sold at a regular price or at a markdown price. The regular price  $r$ is known to the firm from the beginning (i.e., at $t = 0$) and remains constant throughout the planning horizon. In contrast, the retailer must determine the price of reduction $p_{t}$ for the markdown inventory at each time $t$. The price of reduction at time $t+1$ does not exceed the price of reduction at time $t$, which is imposed by the linear constraint 
\begin{equation} \label{con2}
p_{t+1} \leq p_{t}. 
\end{equation}

{\bf The SDP Model:} 
We consider a joint decision-making framework in which the firm simultaneously addresses demand uncertainty and internal competition through coordinated ordering and markdown pricing strategies. Within this framework, the firm must dynamically determine how much fresh product to order, when to reduce the old inventory, and what price to offer for the aged inventory in response to uncertain demand conditions. 
The proposed stochastic dynamic programming (SDP) model determines these decisions (of ordering quantity and markdown prices) over a finite planning horizon to measure trade-off between two conflicting objectives by maximizing the total expected discounted profit and minimizing the expected total cost of ordering, holding, and penalties for shortage and wastage of expired products. 

Given predefined selling prices for both regular and markdown inventories, we analyze two distinct cases to compute expected revenue, each based on different demand realizations. When demand is lower than the available inventory, that is $\displaystyle  \tilde{d_{t}} \leq X_{t} $, the total demand is fulfilled. The revenue obtained from the sale at the regular price at time $t$ becomes $r \tilde{d_{t}}$.
In case of high demand $ \displaystyle  \tilde{d_{t}} \geq X_t$, the revenue is  $  \displaystyle r X_t $. 
Combining these cases, the revenue from the sale of products at the regular price in period $t$ can be calculated as
$ R_t =\displaystyle \mathbb{E}_d \left[ r \, min\left\{  X_t,\;  \tilde{d_{t}} \right\} \right].$
Similarly, when the inventory with age $i$ at time $t$ is reduced to price $ p_{t}$, the firm's revenue is  
$ M_t = p_{t} \mathbb{E}_d \left[ \; min\Big\{  Y_t, \; \tilde{d}_{t}  \Big\} \right] $. Note that $\mathbb{E}_d$ represents the expectation operator to be taken over uncertain customer demand. 
Let $c, h,$ and $\gamma$ represent the unit cost of ordering, holding, and shortage, respectively. The operating cost at time $t$ is calculated as $O_t= Q_t + H_t + S_t $, consisting of the 
ordering cost $Q_t = c q_{t}$, the holding cost 
$ H_t =\displaystyle h \mathbb{E}_d \left[ \Big( X_t + Y_t - \tilde{d_{t}} \Big)^+ \right]$ and the shortage cost
$ S_t =\displaystyle h \mathbb{E}_d \left[ \Big( \tilde{d_{t}} - X_t - Y_t  \Big)^+ \right]$. Here, the maximum function $(a)^+ = max\{ a, 0 \}$ takes a value of $a$ if and only if $a > 0$; otherwise, it is zero. Moreover, we take into account the cost of discarding remaining inventory (so-called wastage penalty) aged $n$ as $ W_t =\displaystyle \eta \mathbb{E}_d \left[ \Big(  x_{nt} + y_{nt}  - \tilde{d_{t}} \Big)^+ \right]$. 

Considering that the price of the marked down inventory cannot be increased or reassigned to the regular price once it is in the deterioration stage, we use information regarding the markdown price $p_{t-1}$ determined at previous time while formulating the dynamic model at time $t$. 
Let $C_{t} \large((\mathbf{I}_{t},\mathbf{w}_{t}) \;| \; {p_{t-1}}\large)$ and $P_{t} \large((\mathbf{I}_{t},\mathbf{w}_{t}) \;|\; {p_{t-1}}\large)$ denote the value functions at state of $(\mathbf{I}_{t}$, $\mathbf{w}_{t}$) given ${p_{t-1}}$. Then, the expected total cost over the planning horizon can be computed by solving the following dynamic programming problem;    
\begin{eqnarray}
EC_{min}:\;\;\;\; C_{t} \Big((\mathbf{I}_{t},\mathbf{w}_{t}) \,| \, {p_{t-1}}\Big)
&=& \underset{\begin{subarray}{c}
0\leq q_{t}\leq \kappa, \\0\leq p_t\leq {p_{t-1}}   \end{subarray}}
\min O_t + W_t + \mathbb{E}_d \left[ C_{t+1} \Big((\mathbf{I}_{t+1},\mathbf{w}_{t+1}) \,| \, {p_{t}}\Big) \right]. 
\end{eqnarray}
Similarly, the expected profit incurred by the system over the planning horizon in state $({\mathbf{I}}_t, {\mathbf{w}}_t)$ is formulated as follows; 
\begin{eqnarray}
EP_{max}:\;\;\;\; P_{t} \Big((\mathbf{I}_{t},\mathbf{w}_{t}) \,|\, {p_{t-1}}\Big)
&=& \underset{\begin{subarray}{c}
0\leq q_{t}\leq \kappa, \\ 0 \leq p_t \leq {p_{t-1}}  \end{subarray}}
\max R_t + M_t - O_{t} + \mathbb{E}_d \left[ P_{t+1} \Big((\mathbf{I}_{t+1},\mathbf{w}_{t+1})\,|\, {p_{t}}\Big) \right]. 
\end{eqnarray}

Let $V_{t} \Big((\mathbf{I}_{t},\mathbf{w}_{t}) \,| \, {p_{t-1}}\Big)$ represent the value function in state $(\mathbf{I}_{t}$, $\mathbf{w}_{t}$) of the system at time $t$, given the markdown price ${p_{t-1}}$ of deteriorating inventories at time $t-1$. We also introduce a parameter $0 < \alpha < 1$ to quantify the relative importance of these two conflicting objectives: maximizing profit and minimizing total cost. Then the value function for the dynamic joint ordering-markdown pricing problem at time $t$ can be formulated as follows; 
\begin{eqnarray}
\begin{aligned}
&EWT_{\alpha}:\; & \\
&  V_{t} \Big((\mathbf{I}_{t},\mathbf{w}_{t}) \,|\, {p_{t-1}}\Big)
=\; \max \;\; \alpha \Big( R_t + M_t - O_{t} \Big)  - 
(1-\alpha) \Big( O_{t} + W_t \Big) +  
\mathbb{E}_d \left[ \delta V_{t+1} \Big((\mathbf{I}_{t+1},\mathbf{w}_{t+1})\,|\, {p_{t}}\Big) \right] & \\ 
& \hspace{3.75cm} \text{s.t.} \; \; \qquad 0\leq q_{t}\leq \kappa, & \\
& \hspace{3.75cm}  \qquad \qquad  0 \leq p_{t} \leq p_{t-1}. \quad &
 \end{aligned}
\end{eqnarray}
where $\delta \in (0, 1]$ is the discount factor for one-period. Since no product is to be sold after the planning horizon, the boundary condition at the terminal time period is expressed in terms of the total holding and disposal costs for the regular and markdown inventories as follows;
\begin{eqnarray}
V_{\tau+1} \Big( (\mathbf{I}_{\tau+1}, \mathbf{w}_{\tau+1} ) \,|\, {p}_{\tau} \Big) &=& - (h+\eta) \sum_{i=0}^{n}(x_{i,\tau+1}+y_{i,\tau+1}). 
\end{eqnarray}
The next section analyzes the structural properties of the SDP model and utilizes these properties for state and action space reduction strategies that allow us to solve efficiently the stochastic dynamic model to derive the optimal policy. 

\section{State and Action Space Reduction Strategies}

 The SDP formulation of the joint ordering and markdown pricing problem exhibits several structural properties that are utilized to develop a solution algorithm. In this section, we  first establish the structural properties of the SDP model based on $k$-concavity and sub-modularity of the value function. Then we introduce state and action space reduction rules using these properties. Theoretical analysis allows us to develop an efficient solution strategy by reducing both the state and action spaces. 
 The single-period profit function $g_{t}\Big(\mathbf{(I_{t}}, \mathbf{w_{t}})\, | \, {{p_{t-1}}} \Big)$ at time $t$ is formulated as follows; 
\vspace{0.5mm}
\begin{equation*}
\begin{aligned}     
g_{t}\Big(\mathbf{(I_{t}},\mathbf{w_{t}})\,|\,{p_{t-1}} \Big)  = &\;  \alpha \Big(R_t + M_t - O_t\Big) - (1-\alpha) \Big(O_t + W_t\Big) \\
=& \; \alpha \mathbb{E}_d \Big[  r\; min\Big \{X_t, \; \tilde{d}_{t} \Big\}  +  p_{t} min \Big\{Y_t, \; \tilde{d}_{t}\Big\} \Big] - \gamma \mathbb{E}_d \Big[max\{\tilde{d}_{t} - X_t - Y_t, 0\} \Big] - c q_{t} \\
& -h (X_{t+1}+Y_{t+1})-(1-\alpha) \eta (x_{n,t+1}+y_{n,t+1}). \nonumber
\end{aligned}
\end{equation*}

\noindent As $k$-concavity of the value function is derived from the concavity of the single-period profit function, we can rewrite the SDP model in terms of $g_{t}\Big(\mathbf{(I_{t}}, \mathbf{w_{t}})\, | \, {p_{t-1}} \Big)$ as
\vspace{0.25cm}
\begin{equation} \label{eqn:g_val}
\begin{aligned}
&  V_{t}\Big(\mathbf{(I_{t}}, \mathbf{w_{t}})\, | \, p_{t-1} \Big) = \underset{\begin{subarray}{c}
0\leq q_{t}\leq\kappa, \\ 0 \leq p_{t} \leq p_{t-1}  \end{subarray}}
{\text{max}} &  g_{t}\Big(\mathbf{(I_{t}}, \mathbf{w_{t}})\, | \, 
p_{t-1} \Big) +  \mathbb{E}_d \Big[ \delta V_{t+1}\Big(\mathbf{ (I_{t+1}}, \mathbf{w_{t+1}}) \, | \, p_{t} \Big) \Big].  
\end{aligned}
\end{equation}

\noindent {\color{black}The following propositions establish the main analytical properties of the model. Propositions 1 and 2 state the concavity of the single-period profit function with respect to the types (regular and markdown) of inventories and ordering decisions, as well as the $k$-concavity of the value function of the SDP model, respectively. All proofs of the propositions and theorems presented in this section are provided in the electronic companion.}

\noindent \textbf{Proposition 1:} \textit{Given the price $p_{t-1}$ of the perishable product at time $t-1$ and a fixed markdown policy $\mathbf{w_{t}}$ at time $t$, the single-period profit function $g_{t}\Big(\mathbf{(I_{t}}, \mathbf{w_{t}})\, | \, p_{t-1} \Big)$
is concave in the inventory level $\mathbf{I_{t}}$ and the amount of order $q_{t}$ at time $t$.} 

\noindent \textbf{Proposition 2:} \textit{Given the price $p_{t-1}$ of the perishable product at time $t-1$ and a fixed markdown policy $\mathbf{w_{t}}$ at time $t$, the value function $V_{t}\Big(\mathbf{(I_{t}},\mathbf{w_{t}} )\,|\, p_{t-1} \Big)$ is $k$-concave in the inventory level $\mathbf{I_{t}}$ and the order quantity $q_{t}$ at time $t$.} 

\noindent \textbf{Proposition 3:} \textit{
Given the price $\mathbf{p_{t-1}}$ of the perishable product at time $t-1$ and a fixed markdown policy $\mathbf{w_{t}}$ at time $t$, the single-period profit function $g_{t}\Big(\mathbf{(I_{t}}, \mathbf{w_{t}})\, | \, p_{t-1} \Big)$ is sub-modular in the inventory level $\mathbf{I_{t}}$ for each age group and in the quantity of orders $q_{t}$ at time $t$.}

\noindent To establish the sub-modularity of the single-period profit function as well as the value function, consider the regular inventory level $x_{kt}=s$ with age  $k$ at time $t$. The corresponding inventory vector is represented as $\mathbf{I}_{t}^{(x_{kt} = s)}=(x_{1t},x_{2t},\cdots,x_{kt}=s,\cdots,x_{nt}, y_{1t},y_{2t},\cdots,y_{kt},\cdots,y_{nt})$.

\noindent \textbf{Proposition 4:} \textit{
Given the price $p_{t-1}$ of the perishable product at time $t-1$ and a fixed markdown policy $\mathbf{w_{t}}$ at time $t$, the value function $V_{t}\Big( \mathbf{(I_{t}}^{(x_{kt})},\mathbf{w_{t}}) \;|\; p_{t-1} \Big)$ is sub-modular in the inventory level $\mathbf{I_{t}}$ and the quantity of orders $q_{t}$ at time $t$, given that the future value function satisfies the following monotonicity condition with respect to the inventory level $x_{k+1,t+1}$; 
\vspace{0.25cm}
\begin{equation*}
    V_{t+1} \Big( (\mathbf{I}_{t+1}^{(x_{k+1,t+1} = 0)}, \mathbf{w_{t+1}})\,|\, p_{t} \Big) \; \leq \;  V_{t+1} \Big( (\mathbf{I}_{t+1}^{(x_{k+1,t+1} = 1)}, \mathbf{w_{t+1}}) \, |\, p_{t} \Big) \; \leq \; \\
     \cdots \leq V_{t+1} \Big((\mathbf{I}_{t+1}^{(x_{k+1,t+1} = s)}, \mathbf{w_{t+1}})\, | \, p_{t} \Big).
\end{equation*} 
}

Note that although we consider only regular inventory in Proposition 1, its result can also be extended to markdown inventories. From Propositions 1 and 2, the summation of the single-period function and the future value functions in (\ref{eqn:g_val}) is strictly concave in ordering decisions. Thus, we can leverage the concavity property of the SDP model in (\ref{eqn:g_val}) to efficiently expedite the search for the optimal quantity of order. This property is stated in the following theorem. 
\vspace{-0.2cm}
\begin{theorem}{\bf (Action Space Reduction):}
Given the price $p_{t-1}$ of the perishable product at time $t-1$,
consider a state $\mathbf{I_{t}^s}$ where the regular inventory level at age $k$ is $x_{kt} = s$ as well as the markdown information $\mathbf{w_{t}}$ related to different aged products at time $t$. 
\begin{compactitem}
\item[(i)] If the maximum value function is achieved at the ordering quantity $q_{t}^{*}=a$ units, then the value function with a higher ordering quantity $q_{t}^{}=e$, where $e \geq a+1$, is always smaller than the one at the optimal order quantity $q_{t}^{*}=a$. Thus, the following relation is valid
\vspace{0.25cm}
$$ V_{t}^{e}\Big(\mathbf{(I_{t}^{s}},\mathbf{w_{t}})\,|\,p_{t-1}\Big) < V_{t}^{a}\Big(\mathbf{ (I_{t}^{s}},\mathbf{w_{t}})\,|\,p_{t-1}\Big),$$
where $V_{t}^{e}\Big(\mathbf{(I_{t}^{s}},\mathbf{w_{t}})\,|\,p_{t-1}\Big)$ and $V_{t}^{a}\Big(\mathbf{(I_{t}^{s}},\mathbf{w_{t}})\,|\,p_{t-1} \Big)$ denote the value functions with the ordering quantities $e$ and $a$, respectively.
 \item[(ii)] Suppose that the optimal ordering quantity is obtained in state $x_{kt} = s$ as $q^{*}_{t}=a$. For inventory levels $x_{kt} = \acute{s}$ for $ \acute{s} \geq s+1$, the value function with an ordering quantity greater than $a$ units will always be smaller than the one obtained at the optimal order quantity. 
\end{compactitem}
\end{theorem} 
\vspace{-0.2cm}
This theorem shows that the action space for searching  for the order quantity can be extended. The value function in each state does not need to be evaluated to order quantities $e\geq a+1$, unlike the conventional method of backward dynamic programming where all potential values are explored exhaustively. Note that the action space for searching the ordering decisions can be further reduced by implementing sub-modularity properties established in Propositions 3 and 4. 

\begin{figure}[h!]
\centering
\includegraphics[scale=0.6]{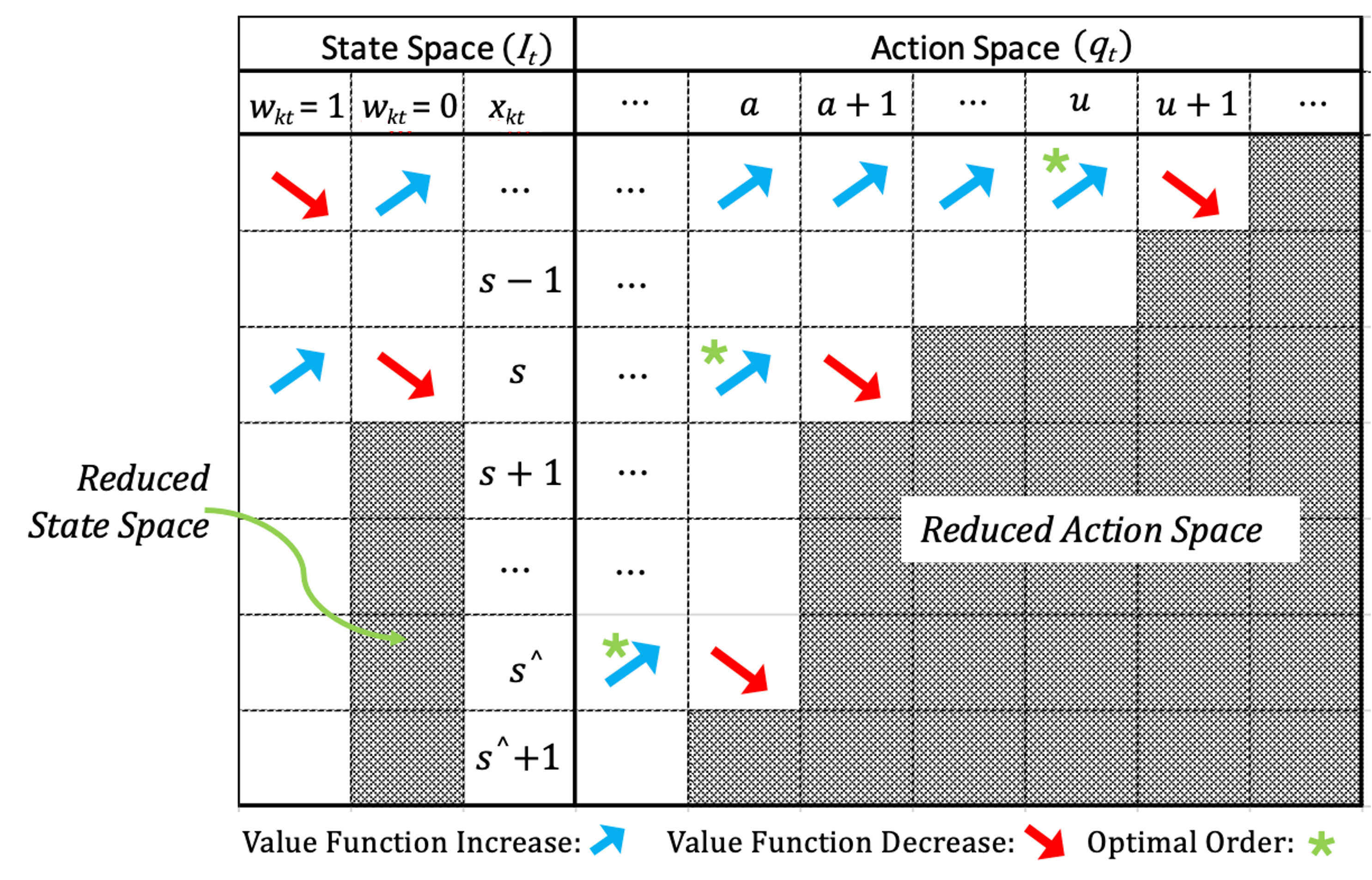}
\caption{Reduction observed at state and action spaces } 
\label{red}
\end{figure}

Given that the value function is sub-modular for different aged inventory levels as well as the order quantity, any increase in the inventory level (while keeping the order quantity constant) would eventually yield a lower value in profit. As stated in part (ii) of Theorem 1, 
the value function is not computed at order quantities higher than the optimal one (that is $e = a+1, a+2, ...., \kappa$) since the relationship
$V_{t}^{e}\Big(\mathbf{(I_{t}^{\acute{s}}},\mathbf{w_{t}})\,|\,p_{t-1} \Big) < V_{t}^{a}\Big(\mathbf{ (I_{t}^{\acute{s}}},\mathbf{w_{t}})\,|\,p_{t-1} \Big),$ holds for all states $x_{kt}>s$. 
Thus, the action search space with respect to the order quantity is progressively trimmed down. Figure \ref{red} illustrates the potential reduction observed in state and action spaces. Notice that parts (i) and (ii) of Theorem 1 correspond to row-wise and column-wise reductions, respectively, in the action space of the ordering quantity.

Next, we focus on state-space reduction through markdown decisions
by exploring the direct relationship between the markdown structure and inventory levels. 
Let us introduce a markdown state variable $\mathbf{w}_{t}^{1}=\Big(w_{1t},w_{2t},\cdots,w_{kt}=1, w_{k+1,t}, \cdots,w_{nt} \Big)$  
for regular inventory $x_{kt}$ with an age of $k$ at time $t$. Similarly, $\mathbf{w}_{t}^{0}=\Big(w_{1t},w_{2t},\cdots,w_{kt}=0,w_{k+1,t},\cdots,w_{nt} \Big)$ represents a no-markdown state. If conducting a markdown sale for $s$ units of inventory aged $k$ at time $t$ is determined to be profitable, then the markdown state remains the same for inventory levels higher than $s$ in the system. The following theorem summarizes this situation and helps to reduce the markdown state space. 
\begin{theorem}{\bf (State Space Reduction): }{Given the price $p_{t-1}$ of the perishable product in period $t-1$, consider a state defined by $\mathbf{I_{t}^s}$, where the regular inventory at age-$k$ is $x_{kt} = s$, together with discount information $\mathbf{w_{t}}$ of products at different ages at time $t$. Suppose that the value function is more profitable in the markdown state ($w_{kt}=1$) than in the no-markdown state ($w_{kt}=0$). If  
$ V_{t}\Big( (\mathbf{I}_{t}^{s},\mathbf{w}_{t}^{1})\,|\,p_{t-1} \Big) > V_{t}\Big( (\mathbf{I}_{t}^{s},\mathbf{w}_{t}^{0})\,|\,p_{t-1} \Big)$ holds, then the value function with the markdown state remains profitable even at a regular inventory level that exceeds $s$. Thus, for regular inventory levels $e \geq s+1$, the following relationship holds: 
$$V_{t}\Big( (\mathbf{I}_{t}^{e},\mathbf{w}_{t}^{1})\,|\,p_{t-1} \Big) > 
V_{t}\Big( (\mathbf{I}_{t}^{e},\mathbf{w}_{t}^{0})\,|\,p_{t-1} \Big).$$
}
\end{theorem} 
\noindent This theorem indicates that there is no need to explore the state of the system without markdown ($w_{kt}=0$) at a higher inventory level ($x_{kt}=s+1$) if higher profit is achieved in state ($w_{kt}=1$) with inventory level $x_{kt}=s$. Therefore, the result  defines a monotonic relationship between profitability and markdown adoption, establishing a benchmark markdown scheme at all inventory levels. Specifically, the optimal markdown scheme at any inventory level can be determined by analyzing the inventory level that is one unit lower than the current one, as described in Theorem 2. If the profitability condition defined in Theorem 2 is satisfied, the set of feasible markdown states can be progressively pruned for the specific inventory level. Consequently, the evaluation of no-markdown states becomes unnecessary beyond the threshold level, leading to a substantial reduction in the effective state space.

\section{Computational Experiments }

\subsection{Design of Experiments}

We design a series of computational experiments to meet two research objectives, i) illustrate the performance of the solution approach, developed based on our theoretical results in Section 4, to solve the dynamic ordering-markdown pricing model described in Section 3, and ii) analyze the benefits of implementing a dynamic joint ordering-markdown pricing strategy for the management of perishable products. The numerical experiments, aimed at evaluating the benefits of the dynamic joint ordering-markdown pricing policy, are specifically designed to address the following managerial questions:

\begin{compactitem}
\item What are the key factors influencing retailers' different goals (maximizing profit and minimizing cost) 
in sale of perishable goods via age-differentiated multiple markdowns?
\item How can the optimal ordering and best-timed markdown strategy 
effectively balance trade-off between profitability and product wastage?
\item How does flexibility in timing and price reduction of multiple markdowns enhance the cost, profitability, and wastage? 
\end{compactitem}

In our numerical experiments, we consider the sale of a perishable product with a three-period lifetime over a finite time horizon of 18 discrete periods. This aligns with the practical assumption that no more than three-age categories coexist on the shelf. Although most of the markdown literature considers only two-age groups (new and old), we extend the shelf capacity rather than introducing further age differentiation. This choice enables a more tractable comparison between our k-concavity based approach and traditional backward dynamic programming. However, we would like to note that the proposed model and solution methodology are easily extendable to settings with additional age differentiation, as the structural properties remain valid for a general $n$-period lifetime product.

Our experimental design considers different parameter settings related to prices and customer segments. Table \ref{table:param-ch3} summarizes the list of parameters for the base case.
\begin{table}[ht]
	\begin{center}
		\small 
		\caption{Selection of parameters for the base case}
		\label{table:param-ch3}
		\begin{tabular}{|l|l|}
			\hline \hline
			Parameters (unit) & Set Values \\	[2.5pt]	
			\hline \hline
			Ordering cost & $c=1$ \\
			Holding cost &  $h=0.01$ \\
			  Shortage cost & $\gamma = 1.1$ \\
             Wastage cost & $\eta= 0.5 $ \\
			Regular price & $r = 4$ \\
			Markdown reductions & 10\%, 15\%, 20\%, 25\% or 30\% \\
			\hline \hline 
		\end{tabular}    
	\end{center}
\end{table}
We adopt the base case from \cite{li2012grocery} and \cite{li2016managing} for our computational study. The ordering, holding and penalty costs in the base case are assumed to be 1, 0.01 and 0.5, respectively. The ordering capacity for the perishable product is assumed to be $\kappa = 15$. The regular price for selling the perishable product is set at 4. Additionally, markdown reductions of  10\%, 15\%, 20\%, 25\% or 30\% are considered for the regular sale price of perishable products 
\textcolor{black}{\citep{gaurd1,wasteless_1}}. 
 We now briefly explain the customer choice model that will be integrated into the optimization model and describe the main steps of the simulation experiments. 

{\bf The Customer Choice Model:}  Price reduction can attract price-sensitive consumers who would not purchase the product at its regular price, but can also induce demand cannibalization by diverting customers who would otherwise buy the product at full price \citep{ferguson2007should,li2012grocery}. To capture this trade-off, we adopt a discrete-choice framework that models how customers select among available purchase options based on their perceived utilities. This model accounts for various observable and unobservable factors that influence a customer's decision to purchase a perishable product \citep{train2009discrete}. Observable factors include the price, quality, and freshness of the product, as well as the customer’s sensitivity to freshness. Unobservable factors, on the other hand, relate to individual preferences such as personal taste and budget. In this study, we consider both observable and unobservable factors to formulate the customer choice model for the joint ordering and markdown pricing problem. 

At each time $t$, a customer faces three alternatives: purchasing the product at the regular price, purchasing it at the markdown price, or leaving the system without purchasing (represented by $s,m,o$, respectively). Each customer derives utility from a product based on its quality (or freshness), price (regular or discounted price), and random preference factors. The inventory available at time $t$ can be divided into two groups: regular and markdown. For ease of notation, we denote the price of option $j\in \{s,m\}$ by $\rho_{jt}$. Although the regular inventory is sold at the fixed regular price $\rho_{st} = r$, $\rho_{mt}$ represents the discount price offered at time $t$. Let $\mathcal{K}$ denote the set of customers. For each customer $k \in \mathcal{K}$, let $U_{jt}^{k}$ denote the random utility associated with choosing the option $j \in \{s,m,o\} $ at time $t = 1,\cdots, \tau$. The utility of the no-purchase option is given as $U_{ot}^{k} = 0$.  Following \cite{train2009discrete} and \cite{ akccay2010joint}, the utility of customer $k$ to choose option $j\in \{s,m\} $ can be expressed as a linear function of price and quality level as follows;  
\begin{equation}
    U_{jt}^{k} =  \theta^{k} \alpha_{jt} - \rho_{jt}   + \mu   \epsilon_{j}^{k} 
\label{eq:utility}
\end{equation}
\noindent where $\theta^{k}$ measures the customer’s sensitivity to quality, $\alpha_{jt}$ denotes the level of quality of the product, the random term $\epsilon_{j}^{k}$ captures idiosyncratic preferences of customer $k$ for choice $j$, and $\mu$ measures the degree of such preference heterogeneity. Assuming that $\epsilon_{j}^{k}$ follows a Gumbel distribution, the probability that customer $k$ selects option $j$ at time $t$ is given by
\begin{equation}
\gamma_{t}^{k}(j) = \frac{ e^{((\theta^{k} \alpha_{jt} - \rho_{jt}  )/\mu) }}{1+ \displaystyle \sum_{j \in \{s,m \}} e^{((\theta^{k} \alpha_{jt} - \rho_{jt}  )/\mu) }  }. 
\end{equation}
Although the quality level and market price of the product are the same for all customers, heterogeneity in $\theta^{k}$ and $\epsilon_{j}^{k}$ leads to individual differences in the probability of purchase. These choice probabilities are subsequently used to determine the expected demand in each period of the stochastic dynamic program.

At each time period, customers encounter two groups of products on the shelf, those sold at the regular price and those sold at the reduced price. Within each group, all products share the same price, which is determined by the retailer’s markdown and ordering decisions. Customers therefore do not distinguish among individual items
within a group (e.g., among items of slightly different ages) but instead evaluate the utility of the group representative based on its common quality level and price.
The utilities derived from \eqref{eq:utility} are also used to simulate customer purchasing behavior over time. In each simulation run, customer-specific parameters such as quality sensitivity $\theta^k$, scale parameter $\mu$, and random preference terms $\epsilon_j^k$ are first initialized for all $k \in \mathcal{K}$. Then, at each time period $t$, the following steps are performed:
\begin{compactitem}
    \item Customer arrivals are generated from a Poisson distribution with a mean rate of $15$ customers per period, and the maximum number of customers in any period is capped at $30$.
    \item For each arriving customer and option $j \in \{s,m,o\}$, utilities are calculated using equation (\ref{eq:utility}) based on the specified prices and quality levels. 
    \item Each customer $k$ chooses the option that provides the highest utility.
    \item Individual choices are aggregated to obtain the realized demand for each option. The resulting demand is then used to update the remaining regular and markdown inventories before proceeding to the next period. 
\end{compactitem}
This process is repeated in all time periods to generate simulated demand trajectories under different pricing and freshness conditions.
Each customer's request depends on their sensitivity to freshness. Thus, we assume that each customer $k$ is assigned a value of $\theta^k \in [0, 1]$, representing their sensitivity to freshness (or quality). A customer with a higher value of $\theta^k$, closer to 1, is highly sensitive to freshness and is classified as sensitive to quality. On the other hand, a customer with a lower $\theta^k$ value is price sensitive. To ensure an equal number of customers in each segment, we uniformly distribute the value of $\theta^k$ in the initial experiments.  We also conduct experiments in which we vary the proportion of customer segments. The expected profits are estimated by simulating customer requests over \textcolor{black}{1000} sample paths.

\subsection{Computational Results }

This section presents numerical results and summarizes the key findings in three categories to address the main operational and managerial questions outlined above. 

{\bf Performance of Solution Approaches:} In the first experiment, the solution approach developed based on our theoretical results in Section 4 is compared with the standard backward dynamic programming (BDP) algorithm.  The proposed solution approach directly applies the rules derived from exploiting the property of $k$-concavity in the dynamic programming model, abbreviated as KCDP. The two methods provide identical optimal profits and policies, confirming the validity of KCDP as an exact solution approach. However, they differ substantially in computational efficiency. Table \ref{table:var-ch3} reports the expected profit, CPU time and the number of evaluations of the state-action value function for both algorithms across different ordering capacities, along with the percentage reduction in evaluations achieved by KCDP.

\begin{table}[ht]
{\footnotesize 
\begin{center}
\caption{Performance comparison solution algorithms } 
\label{table:var-ch3}
\begin{tabular}{|c||c||c|c||c|c||c|}
\hline 
& Optimal &\multicolumn{2}{c||}{\bf BDP Algorithm } & \multicolumn{2}{c||}{\bf KCDP Algorithm}  & Reduction in \\  \cline{3-6}  \cline{3-6}
Capacity&	Expected	& CPU time & Number of  & CPU time& Number of & Value Function   \\ 
(units) & Profit &   (in seconds) & Evaluations & (in seconds) & Evaluations &  Evaluations (\%) \\  [3.5pt] \hline 	\hline		
5  & 188.11 & 5.95      & 126132     & 4.37      & 48096     & 61.87 \\
10 & 407.58 & 323.44    & 857417     & 35.15     & 182147    & 78.76 \\
15 & 629.30 & 3385.81   & 2733952    & 256.37    & 410284    & 84.99 \\
20 & 851.57 & 28241.81  & 6295737    & 2097.05   & 734966    & 88.33 \\
25 & 1074.16& 45528.51  & 12082772 & 6856.71   & 1161457 & 90.39 \\
30 & 1297.07& {2.5 days}  & 20635057        & 18962.45  & 1697108 & 91.78\\
35 & 1520.31 & $NA $ & $NA $ & 21041.76 & 2341919 & $NA $ \\
\hline 	
\end{tabular}    
\end{center}
}
\end{table}

As the ordering capacity increases, both algorithms exhibit longer computation times as a result of the exponential growth of the state and action spaces. However, KCDP consistently achieves significantly lower computational times compared with BDP. To better understand this significant difference in computational times, we examine the state and action space characteristics of both algorithms. The state and action space values represent the number of value function evaluations performed during the solution process where each state requires evaluation over all feasible actions. Consequently, a reduction in evaluations of the value function  directly translates to a decrease in overall computational time. When the ordering capacity is set to 5, KCDP reduces the number of evaluations of the value-function by approximately 60\% relative to BDP. This reduction becomes even more prominent, reaching nearly 80 to 90\%, at higher capacity levels (e.g. 20, 25, and 30). This indicates that KCDP explores substantially fewer state–action combinations leveraging the structural properties of $k$-concavity to eliminate redundant regions of the search space. Therefore, we apply the KCDP algorithm to solve the dynamic programming model in the rest of the numerical experiments. 

{\bf Impact of the Optimal Policy:} In this experiment, we intend to analyze the impact of the optimal policy and derive the key factors (especially the costs of wastage and shortage) that affect two conflicting goals in the management of a perishable product: profit maximization and waste reduction. A firm focused solely on profit maximization seeks to enhance sales; however, strategies aimed at increasing sales often lead to increased wastage due to the perishable nature of the product. Thus, determining the best strategy to balance both objectives is a complex decision-making challenge. In order to illustrate the trade-off between  profit maximization and waste reduction, we first solve the stochastic dynamic programming models $EC_{min}$ and $EP_{max}$ for the given data set. Note that these models can be derived from $EWT_{\alpha}$ by fixing $\alpha = 0$ and 1, respectively. For an equally weighted trade-off between two objectives, we solve the stochastic dynamic programming model  $EWT_{\alpha}$ for a fixed value of ${\alpha=0.5}$. Table \ref{table:var-ch4} illustrates the performance comparison of optimal joint policies in terms of different metrics in two cases (with and without wastage costs) using the base case setting. 
\begin{table}[ht]
\small
\begin{center}
\caption{Performance of different models using base case with and without wastage cost } 
\label{table:var-ch4}
\begin{tabular}{l||rrr||rrr}
\hline 	
\multicolumn{1}{l||}{  } & \multicolumn{3}{|c||}{ \bf Case 1 (wastage cost) } & \multicolumn{3}{c}{ \bf Case 2 (no wastage cost) } \\ [3.5pt] \hline \hline
\multicolumn{1}{l||}{ Metrics (average)} & \multicolumn{1}{|l}{$EP_{max}$} & $EC_{min}$ & $EWT_{\alpha}$ & $EP_{max}$ & $EC_{min}$ & $EWT_{\alpha}$  \\ [3.5pt] \hline \hline
 Revenue (regular price) & 457.24	&	441.17	&	375.96	&	457.24	&	429.93	& 	355.94	\\
 Revenue (markdown price) & 113.95	&	96.85	&	15.41	&	113.95	&	90.69	& 	32.90	\\ 
Expected Total Cost & 171.69	&	153.13	&	107.40	&	179.58	&	150.30	& 	107.82 \\ 
Expected Profit & 399.50	&	384.90	&	283.97	&	399.50	&	370.32	& 	281.01	\\  [3.5pt] \hline  \hline
\multicolumn{1}{l||}{Number of Products (units)} 
 & \multicolumn{1}{|l}{$EP_{max}$} & $EC_{min}$ & $EWT_{\alpha}$ & $EP_{max}$ & $EC_{min}$ & $EWT_{\alpha}$  \\ [3.5pt] \hline \hline
Average Orders & 168.23	&	149.85	&	103.04	&	168.23	&	143.60	& 	101.53		\\
Average Waste & 15.77	&	8.95	&	4.68	&	15.77	&	7.13	& 	3.73		\\
Average Shortage & 2.09	&	2.06	&	3.11	&	2.09	&	1.93	& 	3.29	\\
\hline 	
\end{tabular}    
\end{center}
\end{table}
As expected, the model $EP_{max}$ makes the highest profit, followed by the models of $EWT_{\alpha=0.5}$ and $EC_{min}$. However, $EP_{max}$ also provides the highest waste due to its higher ordering levels. This pattern is observable from both cases regardless of including the wastage cost and highlights the link between excessive ordering levels for profit maximization and the resulting increase in waste. For a firm that wants to minimize waste, the model $EC_{min}$ provides an alternative optimal strategy. Although this model significantly reduces wastage, it does so at the expense of profit, which experiences a notable decline. In contrast, the equally weighted trade-off model  $EWT_{\alpha=0.5}$ serves as a balanced alternative as it achieves approximately a 50\% reduction in waste compared to the $EP_{max}$ model, while profit decreases by only around 7\%. It should also be noted that all models in both cases generate revenue through regular and markdown sales. This observation suggests that a dynamic ordering and markdown pricing strategy can be effectively used to balance the conflicting objectives of profit maximization and waste reduction. 

We next test whether the reduction in wastage arises from explicitly penalizing waste in the objective function or not. To do this, we compare the results obtained by three optimization models using the base case setting with and without waste costs. Since the $EP_{max}$ model does not take waste cost into account, its results remain unchanged in  cases 1 and 2, and comparisons are only made between $EC_{min}$ and $EWT_{\alpha=0.5}$. 

\begin{figure}[h!]
\centering

\includegraphics[scale=0.45]{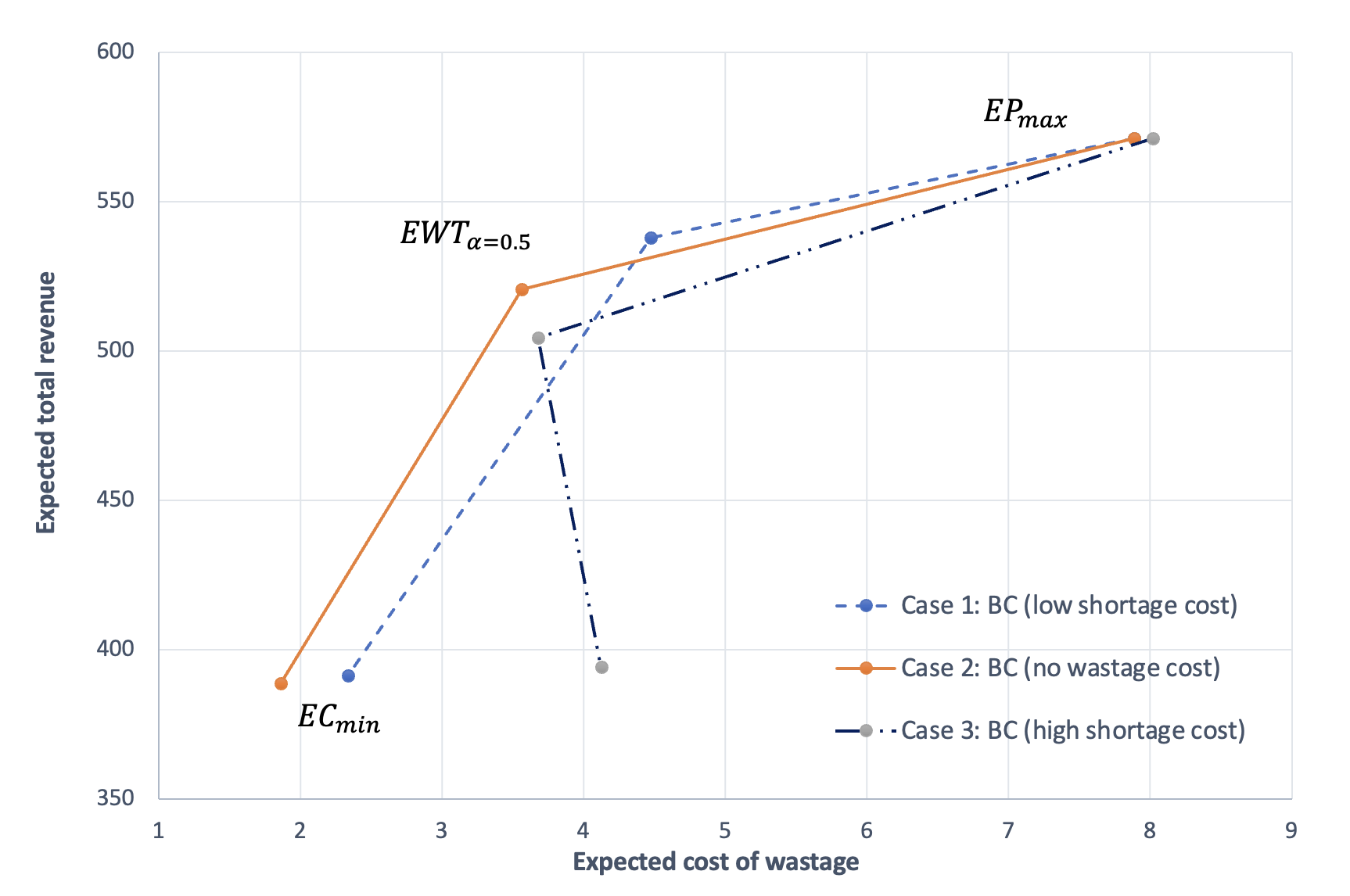}
\caption{Efficient frontiers obtained using the base case setting with different shortage costs }
\label{EF}
\end{figure}

In order to establish impact of the shortage cost on the optimal policy, we increase the unit shortage cost in the base case setting from 1.1 to 1.5 (labeled Case 3). Figure \ref{EF} plots efficient frontiers that display trade-offs in terms of expected revenue and waste using the base case with low and high shortage costs (i.e, cases 1 and 3). 
As can be seen in Figure \ref{EF}, expected profit declines in all models due to higher total costs and a larger quantity of orders intended to minimize unmet demand. Interestingly, 
the model $EC_{min}$ shows a sharp increase in wastage cost, as the elevated shortage cost becomes a dominant component of the total cost, leading to a lower shortage level of just 0.64 units. On the other hand, the minimum wastage levels in case 3 are obtained by $EWT_{\alpha=0.5}$ due to its well-balanced structure between the two conflicting objectives of profit maximization and total cost minimization. 
As another key insight, we can say that the model $EWT_{\alpha=0.5}$ enables the retailer to reduce waste by nearly 50\% while sacrificing only about 5–10\% of the maximum profit attainable for the given test cases. This represents a well-balanced strategy that achieves substantial waste reduction without a significant loss in profitability. 
This approach offers a rare alignment between economic and environmental objectives, supporting sustainable operational decisions. This balance can be achieved through a strategic and dynamic ordering–markdown policy. In the next section, we further explore how markdown strategies operate in all three model configurations. 

{\bf Dynamic Markdown and Ordering Strategies:} Supermarkets often present markdowns as a strategy to reduce product waste, appealing to the sense of sustainability and social responsibility of customers. However, it remains unclear whether this narrative reflects real waste reduction or simply serves as a marketing justification for profit-driven pricing. Indeed, several academic and practitioner-oriented studies suggest that markdowns are primarily motivated by the desire of firms to maximize profits and can, in some cases, exacerbate waste. To address this ambiguity, we design an experiment that specifically examines the role of markdowns in balancing two potentially conflicting objectives: profit maximization and total cost minimization. As discussed above, waste levels vary depending on which of these objectives is prioritized. Based on this insight, we specifically analyze how markdown and ordering strategies differ between three models: $EP_{max}$, $EC_{min}$, and $EWT_{\alpha =0.5}$. 

Industry practices within the retail sector exhibit a wide range of discount strategies \citep{wasteless_2}. Traditionally, markdown policies have been characterized by a single time price reduction implemented shortly before the expiration of the product. However, more recent evidence, as documented by McKinsey \citep{Bahr22}, indicates a strategic shift toward earlier and more frequent price adjustments over the product life-cycle. Empirical findings suggest that even relatively modest early discounts can substantially accelerate sales and mitigate waste. 

In addition, the reports by  \cite{wasteless_3} and \cite{trellis} indicate that Wasteless employs modest early discounts and intensifies markdowns only if products remain unsold, thereby avoiding drastic last-minute reductions. 
Motivated by these evolving industry practices, our experimental setting examines the impact of alternative markdown strategies within a dynamic pricing framework. We consider a retailer managing inventory across three age categories: newly ordered units ($q_t$), 1-period-old inventory ($x_{1t}, y_{1t}$), and 2-period-old inventory ($x_{2t}, y_{2t}$), offered at regular and markdown prices, respectively. Within this setup, the firm can implement different types of multiple reductions as described below:
\begin{compactitem}

\item{\it No Markdown:} 
Products of different ages are offered at regular price, without any price discounts.

\item{\it Complete Markdown:} The pricing policy is characterized by the application of discounts to all age‑differentiated products relative to their regular price benchmarks.

\item{\it Delayed (late) Markdown:} Two‑period‑old inventory that was not discounted in prior periods is now marked down. This action constitutes the last opportunity to avoid waste.

\item {\it Early Markdown:} Discounting is introduced early by reducing the prices of inventories of 1‑period‑old and 2‑period‑old.

\end{compactitem}

\begin{figure}[h!]
\centering
\includegraphics[scale=0.17]{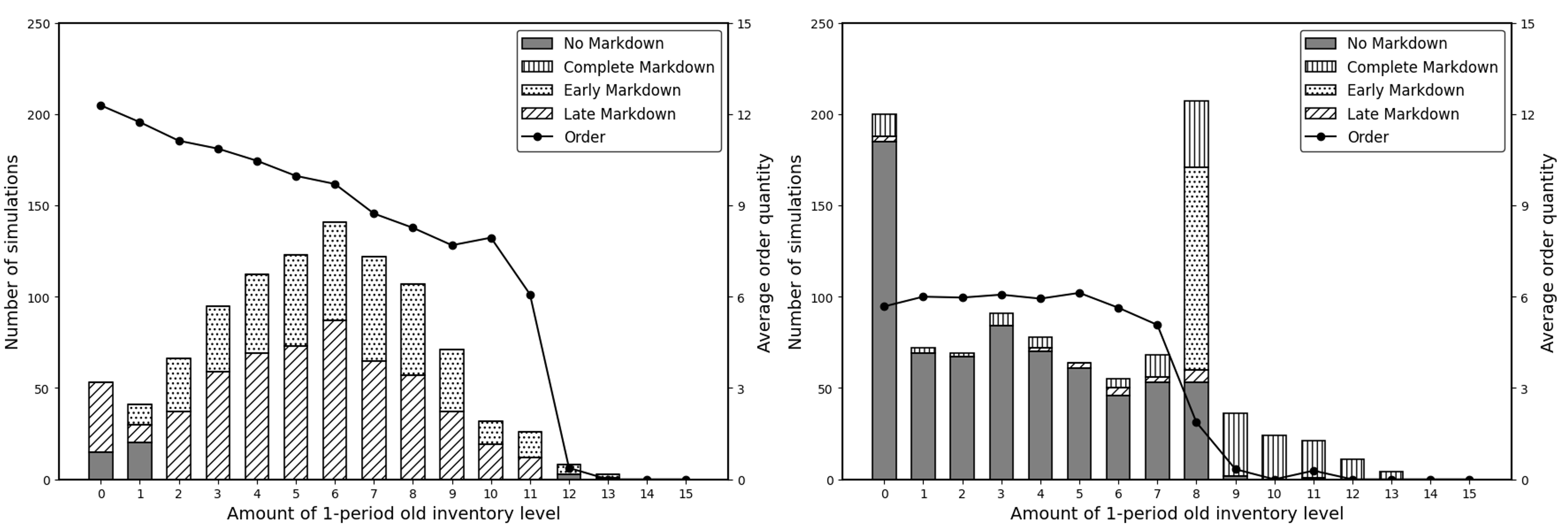}
\caption{Multiple markdown and ordering strategies obtained by $EP_{max}$ (left) and 
$EC_{min}$ (right) }
\label{markdown1}
\end{figure} 

\begin{figure}[h!]
\centering
\includegraphics[scale=0.45]{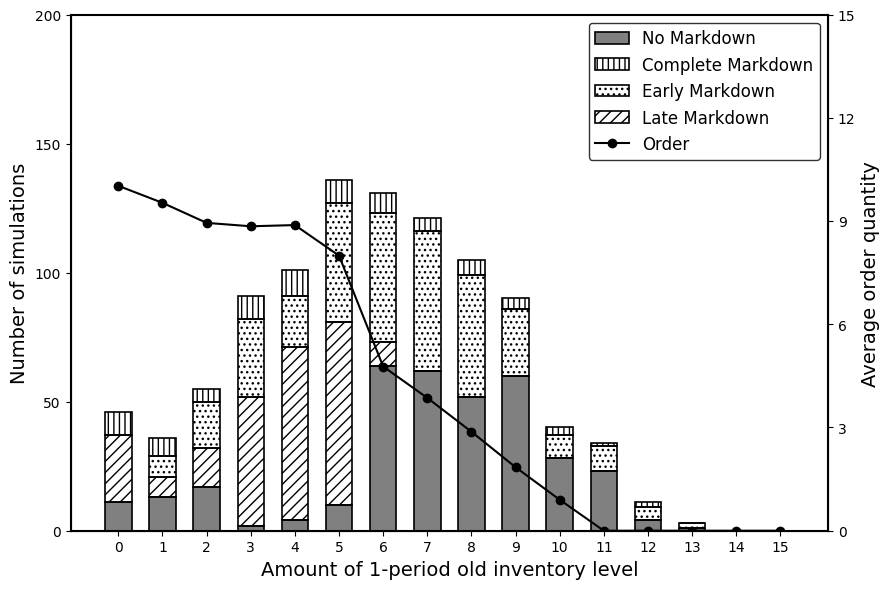}
\caption{Multiple markdown and ordering strategies obtained by $EWT_{\alpha=0.5}$ }
\label{markdown2}
\end{figure}

Figure \ref{markdown1} illustrates different markdowns and ordering strategies derived from two objectives: maximizing expected profit ($EP_{max}$ on the left) and minimizing expected cost ($EC_{min}$ on the right).
Similarly, Figure \ref{markdown2} shows the results of the model $EWT_{\alpha =0.5}$.
These strategies are analyzed in relation to the available inventory of one-period-old units at the final time period, based on 1,000 simulation runs.
The x-axis represents different inventory levels, while the bar charts display the frequency of each markdown strategy at each inventory state over all simulations. 
Each bar is color-coded according to a legend that identifies the specific markdown strategy. The left y-axis indicates the number of simulations (out of 1,000) in which each strategy was selected for a given inventory level. In parallel, the right y-axis shows the average ordering quantity, depicted by a solid-dotted line chart.
For instance, from the graph of $EP_{max}$ shown on the left, we can see that 
an inventory level of six units of stock of one-period-old is observed approximately 140 times, indicating that the system occupied this state 140 times during the simulation. 
The distribution of bars reveals that a wide range of inventory states are visited under the $EP_{max}$ model. 

We observe that partial markdowns dominate in $EP_{max}$ whereas no markdown is relatively rare and  occurs primarily when the inventory of one-period old is very low (zero or one unit) or, occasionally, when it reaches around 12 units. 
The complete markdown  strategy does not occur in this model. These findings suggest that when a firm aims solely to maximize profit, it tends to prefer partial markdown strategies. Such strategies display regular and discounted prices to customers: fresh, high-quality products remain at regular prices, while older inventory is sold at reduced prices. This dual-pricing mechanism enables the firm to effectively cater to different customer segments seeking freshness and quality, as well as those who are more price-sensitive. In terms of ordering, as indicated by the solid line on the secondary y-axis, higher order quantities are observed at lower inventory levels, with ordering gradually decreasing as inventory accumulates.

Next, we examine the $EC_{min}$ model, which demonstrates the lowest frequency of markdowns. Unlike other approaches, no single markdown strategy dominates in this cost-minimization framework. As shown in previous results, $EC_{min}$ consistently achieves lower wastage costs compared to the $EP_{max}$ model in all scenarios.
Interestingly, the observed reduction in waste is achieved primarily through the limited use of markdowns and a strong preference for the no-markdown strategy. This outcome is somewhat counterintuitive, given that supermarkets often rely on markdowns as a tool to minimize waste. However, our findings suggest that avoiding markdowns altogether may be a more effective approach to waste minimization.
To better understand this behavior, we compare the ordering patterns of the $EC_{min}$ and $EP_{max}$ models. Even in a state of zero inventory, $EC_{min}$ places orders of approximately 5 to 6 units—roughly half the quantity ordered under $EP_{max}$. 
Additionally, the reduction in ordering begins earlier in $EC_{min}$, around eight units of one-period-old inventory, whereas in $EP_{max}$ it starts later around ten or eleven units. 
This suggests that $EC_{min}$ follows a more conservative ordering policy, resulting in fewer leftover units that require markdowns. Among the limited reductions observed in $EC_{min}$, the complete reduction strategy is predominant, with the aim of completely 
eliminating the remaining stock. 
Partial markdowns are rarely used. In contrast, $EP_{max}$ deliberately maintains regular and marked-down pricing to appeal to a larger customer base.
In contrast, $EC_{min}$  focuses exclusively on minimizing expected costs through lower order quantities, without offering markdown opportunities preferred by price-sensitive customers. For managers, this implies that while the $EP_{max}$ model can produce an effective profit-oriented strategy, it comes at the expense of greater waste. In contrast, adopting $EC_{min}$  significantly reduces waste, but results in a considerable loss of profit  (as shown in Table 5) and can limit customer engagement due to the absence of markdown pricing. 

Given these trade-offs, a key question arises: can $EWT_{\alpha =0.5}$ serve as a viable middle ground between the two extremes? In order to address this question, we now turn to our analysis of the model to assess whether it can effectively balance the dual objectives of profit maximization and cost minimization through its markdown strategy. As shown in Figure \ref{markdown2}, both the distribution of inventory states and the different types of discount strategies are well balanced. Unlike $EP_{max}$ and $EC_{min}$ (dominated by partial markdowns and no markdown, respectively), 
$EWT_{\alpha =0.5}$ does not have a single dominant policy. Instead, it deploys a healthy blend of all types of markdown strategy, namely partial markdowns and no-markdown, with even occasional occurrences of complete markdown. 

The ordering policy helps explain this mix. The average order quantity in $EWT_{\alpha =0.5}$ exceeds that in $EC_{min}$ but remains below $EP_{max}$.  When inventories are low, the equally weighted trade-off model mirrors $EP_{max}$ in its ordering behavior and mainly adopts 
{\color{black} \it early markdown} by 
selling fresh items at regular price and discounting older units. This enables the firm to serve both quality-seeking and price-sensitive customers. As inventories accumulate (beyond 6–7 units of one-period-old stock), the quantity of orders drops to the same levels as witnessed in $EC_{min}$, and the model changes to 
{\color{black} \it no markdown or early markdown} strategies
to limit excess and accelerate clearance. This diversified portfolio of markdowns and ordering decisions allows the equally weighted trade-off model to simultaneously preserve profit, limit waste, and expand the reach of the  customer. Given its favorable performance, we would like to examine the markdown mechanics of this model further, specifically, how the multiple price reductions unfold within each markdown strategy and how these policies operate in the presence of two-period-old inventory.

\begin{figure}[h!]
\centering
\includegraphics[scale=0.4]{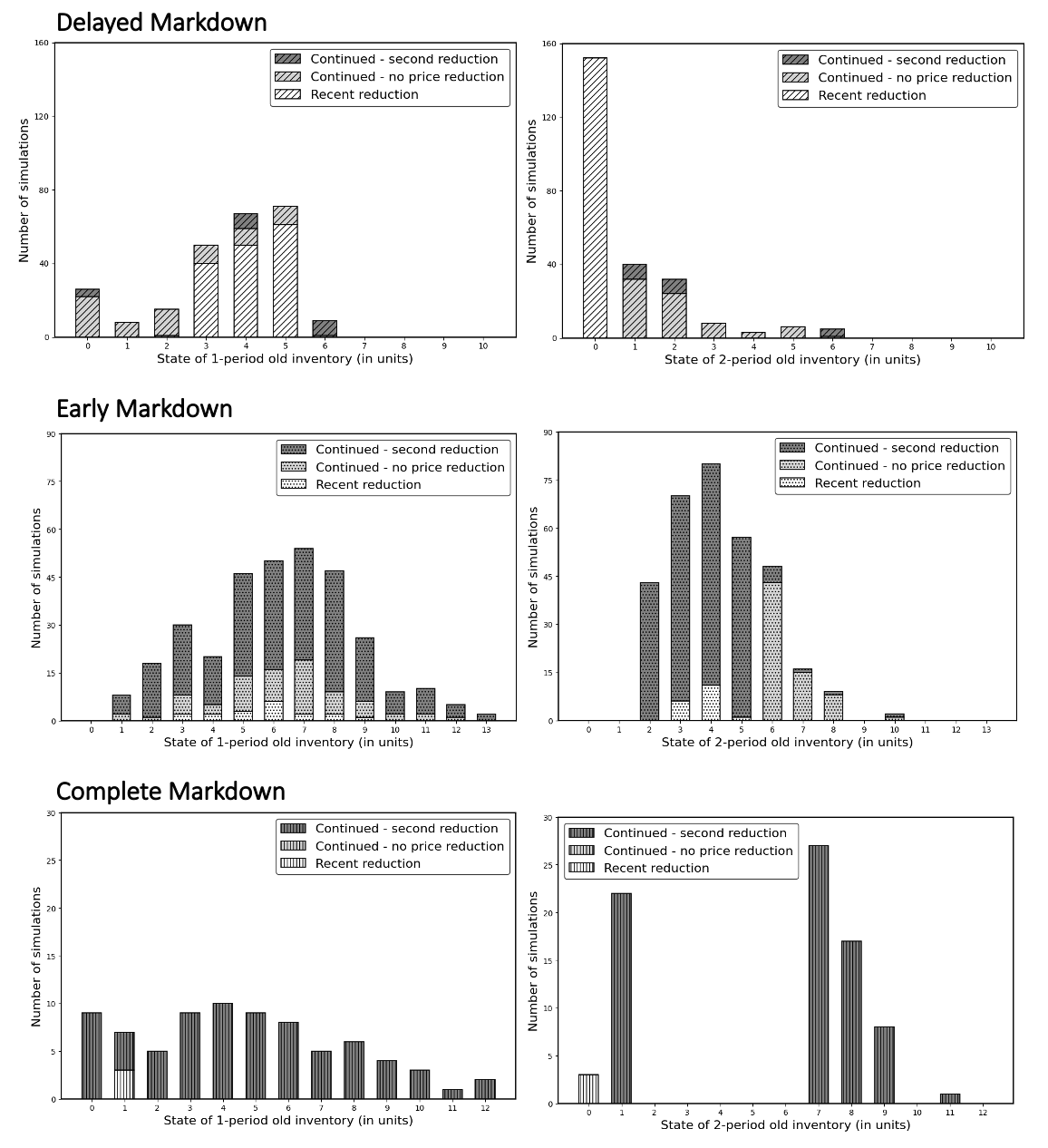}
\caption{ Different variations of multiple markdown policies obtained by $EWT_{\alpha =0.5}$ using  one-period-old (left panels) and two-period-old (right panels) inventories  across all simulation runs}
\label{markdown3}
\end{figure} 

Given its well-diversified portfolio of markdown policies, we analyze 
{\color{black} delayed, early, and complete} markdown strategies
with respect to the states of both one-period-old and two-period-old types of inventories. The resulting relationships are illustrated in six panels in Figure \ref{markdown3}. 
Each panel represents the interaction between inventory age and the corresponding markdown policy. Within each markdown, three variations are highlighted in the bar charts. 
The first variation, called {\it `Recent reduction'}, refers to discounting newly applied to inventory during the current period and is depicted in light gray color all panels. Since the model explicitly tracks markdowns over time, these charts also display the temporal and gradual progression of different markdown policies. If markdowns initiated in previous periods persist into the current one, two additional situations may occur: (i) the markdown continues without additional price reduction, labeled {\it `Continued – no price change'} (shown in medium gray), or (ii) the markdown undergoes a second price reduction, labeled {\it `Continued – second reduction'} (shown in dark gray). 
The relative frequency of each markdown strategy in the overall simulation results 
is displayed in the y-axis. For instance, since {\color{black} \it complete markdown} occurs the least frequently, as observed earlier in Figure \ref{markdown2}, its y-axis is adjusted to reflect its lower incidence.
Recent discounting occurs predominantly under the {\color{black} \it early markdown} strategy 
for inventories of one-period-old (as shown in the  middle row at left panel of Figure 6). 
Note that 
the {\color{black} recent discounting} appears only more frequently when inventory levels from the first period range between 3 and 5 units. The sharp increase in recent markdowns at zero units of two-period-old inventory (as seen in the middle row at the right panel of Figure 6) indicates that when the oldest inventory is absent, early markdowns are applied to the fresher, one-period-old stock. This suggests that, in many cases, the firm aims to prevent the accumulation of two-period-old inventory by initiating markdowns earlier in the product life cycle. However, there are still several instances where two-period-old stock remains unsold. When such an inventory exists, the continued markdowns become more prevalent. 
This implies that markdowns initiated in earlier periods persist due to remaining unsold stock. Among these continued price reductions, the second reduction dominates in both 
{\color{black} \it early and complete markdown strategies}. 
In other words, the firm rarely maintains the same markdown price across periods. Instead, it tends to apply a larger second reduction to expedite sales. This behavior is consistent with practical retail observations and reflects the firm’s intent to accelerate the clearance of older inventory while simultaneously addressing price-sensitive demand and minimizing waste.

Overall, this experiment provides two key insights. First, the equally weighted trade-off model stands out not only because it balances profit maximization and cost minimization, but also because it maintains the well diversified set of markdowns, effectively serving both price-sensitive and quality-oriented customers. Second, the results highlight the importance of allowing multiple price reductions within a dynamic markdown framework. The combination of varying markdown types with multiple price reductions shows that this is an effective mechanism for managing perishable inventory efficiently. 

\section{Conclusions}

Retailers are increasingly deploying computational tools, such as electronic shelf labels and mobile applications, to manage markdowns based on product expiration dates and real-time inventory conditions (Hansen et al., 2024). These emerging technologies and retail innovations reflect the growing shift away from static pricing toward more responsive, data-driven retail practices. As dynamic pricing technologies continue to gain traction in retail, there is a growing need for advanced modeling approaches that extend beyond traditional simple, one-time markdowns toward more flexible, data-driven pricing strategies. These developments indeed underscore the practical relevance of our proposed framework, which enables age-based pricing that adjusts dynamically as products deteriorate. Motivated by this evolving retail landscape, we consider a dynamic markdown management approach that explicitly incorporates the age of the product and the levels of inventory, providing a flexible and efficient solution to the increasing complexity of contemporary retail operations. 

We develop a stochastic dynamic programming model to support retailers in determining optimal ordering and markdown policies for perishable products under demand uncertainty. In contrast to traditional approaches that apply a single markdown close to the expiration date, our model explicitly accounts for the evolution of prices as products age. It allows for multiple age-dependent markdowns and dynamically updates ordering decisions as inventory levels change. This structure enables managers to respond proactively to slow-moving inventory, rather than resorting to late, aggressive price cuts that often result in avoidable waste.  Our formulation therefore generalizes the widely observed two-price structure (regular price and markdown) by allowing retailers to determine both the timing and magnitude of markdowns. This flexibility supports the systematic evaluation of dynamic policies, such as initiating markdowns when inventory reaches a certain age and dynamically replenishing the stock. 
To address the computational challenges inherent in such rich decision environments, we further introduce a rule-based reduction of the state and action spaces. This approach yields optimal dynamic markdown pricing and ordering policies while explicitly capturing the trade-off between revenue maximization and waste reduction. In general, the proposed framework provides analytical insight into optimal pricing policies and offers a timely contribution to the evolving literature and practice of retail markdown management in an era of increasing technological sophistication.

From a managerial perspective, the model highlights an important trade-off between profitability and waste reduction. We explicitly evaluate three strategic regimes. A profit-driven strategy prioritizes revenue maximization by encouraging aggressive markdowns and larger replenishment quantities. Although this approach increases sales, it also increases the risk of overordering and higher levels of expired products. In contrast, a waste-minimization strategy significantly reduces spoilage by limiting inventory exposure, but does so at the cost of foregone sales and lower profits.
In addition to determining optimal dynamic markdown strategies for perishable products, we explicitly account for two competing managerial objectives: maximizing expected profit and minimizing expected total cost, including waste-related losses. These conflicting objectives are incorporated through a weighted objective function, which allows us to examine different strategic priorities. 

Under the profit-maximization regime (EPmax), the model employs aggressive markdown strategies to capture as much demand as possible. While this approach increases expected revenue, it also encourages higher order quantities, which can result in substantial product waste. In contrast, the cost-minimization regime (ECmin) prioritizes waste reduction, significantly reducing spoilage, but at the expense of considerable profit loss. 
The most compelling insight for managers comes from the balanced strategy, which places equal weight on profit and waste considerations (EWT). Our numerical results show that this approach can reduce product waste by nearly 50\% while sacrificing only 5–10\% of the maximum profit achievable. 
This trade-off demonstrates that dynamic markdown and ordering policies can be designed to achieve a win–win outcome, where meaningful waste reduction is achieved while maintaining strong economic performance. This finding underscores a rare alignment between economic performance and environmental objectives, supporting more responsible and sustainable operational decision-making.
The numerical results further emphasize the importance of allowing multiple price reductions within a dynamic markdown framework. Combining various types of markdown with flexible pricing across product age classes proves to be an effective mechanism for managing perishable inventory efficiently. In particular, the equally weighted trade-off model maintains a diverse set of markdown strategies, enabling retailers to simultaneously pursue profit maximization and cost minimization while effectively serving both price-sensitive and quality-conscious customer segments.

Overall, we can conclude that the framework provides managers with actionable guidance on how to use modern pricing technologies more effectively. By combining dynamic markdowns with responsive ordering policies, retailers can improve inventory efficiency, reduce waste, and maintain strong profitability. The results suggest that investments in dynamic pricing capabilities are most valuable when combined with structured decision rules that explicitly balance commercial and sustainability objectives.

\bibliographystyle{agsm}
\bibliography{biblography}

\end{document}